\documentclass[12pt,oneside,reqno]{amsart}
\usepackage{txfonts}
\usepackage{bbm}
\usepackage{amsmath}
\usepackage{graphicx}
\usepackage{mathrsfs}
\usepackage{stmaryrd}
\usepackage{amsfonts}
\usepackage{enumerate,amsmath,amssymb,amsthm}
\numberwithin{equation}{section}

\newcommand{\be}{\begin{eqnarray}}
\newcommand{\ee}{\end{eqnarray}}
\newcommand{\ce}{\begin{eqnarray*}}
\newcommand{\de}{\end{eqnarray*}}
\newtheorem{theorem}{Theorem}[section]
\newtheorem{lemma}[theorem]{Lemma}
\newtheorem{remark}[theorem]{Remark}
\newtheorem{definition}[theorem]{Definition}
\newtheorem{proposition}[theorem]{Proposition}
\newtheorem{Examples}[theorem]{Example}
\newtheorem{corollary}[theorem]{Corollary}

\def\eps{\varepsilon}

\def\a{\alpha}

\def\p{\partial}

\def\[{{\Big[}}
\def\]{{\Big]}}
\def\<{{\langle}}
\def\>{{\rangle}}
\def\({{\Big(}}
\def\){{\Big)}}

\def\bx{{\mathbf{x}}}

\def\dif{{\mathord{{\rm d}}}}

\def\min{{\mathord{{\rm min}}}}

\def\no{\nonumber}
\def\={&\!\!=\!\!&}
\def\bt{\begin{theorem}}
\def\et{\end{theorem}}
\def\bl{\begin{lemma}}
\def\el{\end{lemma}}
\def\br{\begin{remark}}
\def\er{\end{remark}}
\def\bd{\begin{definition}}
\def\ed{\end{definition}}
\def\bp{\begin{proposition}}
\def\ep{\end{proposition}}
\def\bc{\begin{corollary}}
\def\ec{\end{corollary}}
\def\bx{\begin{Examples}}
\def\ex{\end{Examples}}

\def\cB{{\mathcal B}}

\def\cE{{\mathcal E}}
\def\cF{{\mathcal F}}

\def\cL{{\mathcal L}}

\def\cN{{\mathcal N}}

\def\mB{{\mathbb B}}

\def\mE{{\mathbb E}}

\def\mN{{\mathbb N}}

\def\mR{{\mathbb R}}

\def\mT{{\mathbb T}}

\def\sA{{\mathscr A}}

\def\sC{{\mathscr C}}

\def\sG{{\mathscr G}}

\def\sL{{\mathscr L}}

\def\sT{{\mathscr T}}

\def\geq{\geqslant}
\def\leq{\leqslant}

\def\div{\mathord{{\rm div}}}
\def\iint{\int^t_0\!\!\!\int}

\allowdisplaybreaks

\begin{document}

\title{Stochastic Flows of SDEs with Irregular Coefficients and
Stochastic Transport Equations}

\date{}

\author{Xicheng Zhang}

%\thanks{$*$ This work is supported by ARC Discovery grant DP0663153 of Australia.}
%Galerkin's Approximation, Analytic Semigroup,

\dedicatory{
Department of Mathematics,
Huazhong University of Science and Technology\\
Wuhan, Hubei 430074, P.R.China,\\
School of Mathematics and Statistics\\
The University of New South Wales, Sydney, 2052, Australia\\
Email: XichengZhang@gmail.com
 }

\begin{abstract}
In this article we study (possibly degenerate) stochastic differential equations (SDE)
with irregular (or discontiuous) coefficients, and prove that under certain conditions on the coefficients,
there exists a unique almost everywhere stochastic (invertible) flow associated with the SDE
in the sense of Lebesgue measure. In the case of constant diffusions and BV drifts,
we obtain such a result by studying the related stochastic transport equation.
In the case of non-constant diffusions and Sobolev drifts,
we use a direct method. In particular, we extend the recent results on ODEs
with non-smooth vector fields to SDEs. Moreover, we also give a criterion for the
existence of invariant measures for the associated transition semigroup.
\end{abstract}

\maketitle
\rm

\section{Introduction}

Consider the following It\^o's stochastic differential equation (SDE):
\be
\dif X_t=b(X_t)\dif t+\sigma(X_t)\dif W_t,\ \ X_0=x,\label{SSDE}
\ee
where $b:\mR^d\to\mR^d$ and $\sigma:\mR^d\to\mR^d\times\mR^m$ are two
Borel measurable functions, and $(W_t)_{t\geq 0}$
is the $m$-dimensional standard Brownian motion on the classical Wiener
space $(\Omega,\cF,P; (\cF_t)_{t\geq 0})$, i.e., $\Omega$ is the space of all
continuous functions from $\mR_+$ to $\mR^m$ with locally uniform convergence topology,
$\cF$ is the Borel $\sigma$-field, $P$ is the Wiener measure, $(\cF_t)_{t\geq 0}$
is the natural filtration generated by the coordinate process $W_t(\omega)=\omega(t)$.

It is by now a classical result that if $b$ and $\sigma$ are globally Lipschitz continuous, then
there exists a unique bi-continuous solution $(t,x)\mapsto X_t(x)$ to SDE (\ref{SSDE}) such that
for almost all $\omega$ and any $t\geq 0$, $x\mapsto X_t(\omega,x)$ is a homeomorphism.
Thus, $\{X_t(x),x\in\mR^d\}_{t\geq 0}$ forms a stochastic homeomorphism flow (cf. \cite{Ku}).
Recently, there are increasing interests for studying the stochastic homeomorphism flow
property associated with SDE (\ref{SSDE}) under various {\it non-Lipschitz} assumptions on
$b$ and $\sigma$ (cf. \cite{Ma, Ai-Re, Le-Ra, Zh, Fa-Zh,Fa-Zh1, Fa-Lu,Fa-Im-Zh,Zh2}, etc.).
Here, the non-Lipschitz conditions may be less smooth or not global Lipschitz .

On the other hand, when $\sigma$ is non-degenerate and $b$ is not continuous and even singular,
SDE (\ref{SSDE}) may have a unique strong solution for each starting point $x\in\mR^d$
(cf. \cite{Gy-Ma,Kr-Ro, Zh3}, etc.).
But it is not known whether it still defines a stochastic homeomorphism flow.
In the completely degenerate case ($\sigma=0$), a celebrated theory established
by DiPerna and Lions \cite{Di-Li} says that ordinary differential equation (ODE)
\be
\dif X_t=b(X_t)\dif t,\ \ X_0=x\label{ODE1}
\ee
defines a regular Lagrangian flow in the sense of Lebesgue measure when $b$ is a Sobolev vector
field with bounded divergence. This theory was later extended to the case of
BV vector fields by Ambrosio \cite{Am}. The central of DiPerna and Lions' theory are based on
the connection between ODE and the Cauchy problem for the transport equation:
\be
\p_t u+b^i\p_i u=0,\ \ u|_{t=0}=u_0.\label{TR}
\ee
Here and below, we use the usual convention: the repeated indices will be summed.
By introducing a new notion of renormalized solutions, DiPerna and Lions showed
the uniqueness and stability of $L^\infty$-distributional solutions for (\ref{TR}) when $b$ is Sobolev regular so that
they can go back to ODE and show the well posedness of (\ref{ODE1}) with Sobolev vector field $b$
in the distributional sense.

We now back to SDE (\ref{SSDE}). It is also well known that SDE (\ref{SSDE})
is connected with the following stochastic transport equation (cf. \cite{Ku, Ro}):
\be
\dif u=\frac{1}{2}\sigma^{il}\sigma^{jl}\p^2_{ij}u
-(b^i-\sigma^{jl}\p_j\sigma^{il})\p_i u\dif t-\sigma^{il}\p_i u\dif W^l_t
,\ \ u|_{t=0}=u_0.\label{STR}
\ee
Thus, it is natural to ask whether we can extend the DiPerna and Lions theory to the case of SDEs.
Notice that (\ref{STR}) is always a degenerate second order stochastic parabolic equation
whatever $\sigma$ is or not degenerate. More general second order
linear stochastic partial differential equation has been recently studied  in \cite{Zh2}.
In general, it is hard to solve equation (\ref{STR})
if $b$ and $\sigma$ are not smooth (cf. \cite{Ro}).
The source of difficulty clearly comes from the degeneracy. Nevertheless,
we can extend the well known theory about the transport equation
to the case of constant $\sigma$ and BV vector field $b$.
In this case, it will be shown that we can also go back to SDE (\ref{SSDE}) from
stochastic transport equation (\ref{STR}) and obtain the well posedness of SDE (\ref{SSDE})
with BV drift. We remark that in another direction,
Flandoli, Gubinelli and Priola \cite{Fl-Gu-Pr} studied the well posedness of (\ref{STR})
when $b$ is H\"older continuous and $\sigma$ is the unit matrix,
where their proofs  benefit from the stochastic flow associated with SDE (\ref{SSDE}).
We emphasize that when $\sigma$ is constant, SDE (\ref{SSDE}) can be directly solved by
transferring it to a time dependent ODE. But, this will lose some ``stochastic flavor''.

Recently, Crippa and De Lellis \cite{Cr-De-Le} derived some new estimates for ODEs with Sobolev
coefficients. These estimates allowed them to give a direct and simple treatment
for DiPerna-Lions flows. The key ingredient of their method is to give some control for the following
quantity in terms of $\|\nabla b\|_{L^p}$ ($p>1$):
$$
\int_{B_R}\sup_{t\in[0,T]}\sup_{r\in[0,2R]}
\left[\fint_{B_r}\log\left(\frac{|X_t(x)- X_t(x+y)|}{\delta}+1\right)\dif y\right]^p\dif x,
$$
where $B_r:=\{x\in\mR^d: |x|\leq r\}$ denotes the ball with radius $r$ and center $0$.
For estimating this quantity, the Hardy-Littlewood maximal function was used to control the
difference $|b(X_s(x))-b(X_s(x+y))|$.
Moreover, the stability was also derived in \cite{Cr-De-Le} by using a similar quantity.
We remark that the above quantity was first introduced in \cite{Am-Le-Ma} in order to prove the
approximative differentiability of regular Lagrangian flows.
The second part of this paper is to extend Crippa and De Lellis' result to
the stochastic case so that $\sigma$ can be non-constant.

We also mention that Figalli \cite{Fi}  has already developed a stochastic counterpart for
DiPerna-Lions theory. Therein, the martingale solution (or weak solution)
in the sense of Stroock-Varadhan was considered corresponding to the Fokker-Planck equation.
Moreover, the non-degenerate condition on $\sigma$ is required when $\sigma$ is non-constant.
Compared with \cite{Fi}, we can directly construct the ``strong'' solution of SDE (\ref{SSDE})
with Sobolev drift and possibly degenerate diffusion coefficients in the sense of Lebesgue measure.
Moreover, as an easy consequence, we can uniquely solve the SDE in the classical sense when the
initial value is an absolutely continuous $\cF_0$-measurable random variable
(see Corollary \ref{Cor} and Corollary \ref{Cor1} below). It should be noted that for the simplicity,
we only consider the time independent coefficients in the present paper.
Clearly, our results can be extended to the time dependent case by requiring
some integrability in the time variable.

In the study of stochastic dynamical systems, an important problem is to prove
the existence of equilibrium point (invariant measure).
Since we are dealing with non-smooth stochastic differential equations,
it is not expected to have the Feller property for the associated transition semigroup.
Thus, it seems that the classical coercivity condition is not enough to guarantee the existence of
an invariant probability measure for SDE (\ref{SSDE}) (cf. \cite{Ar, Ku}). In the present paper,
we shall give a criterion for the existence of an invariant probability measure in
terms of the classical coercivity condition as well as some divergence condition (see Theorem \ref{Main1} below).
We want to emphasize that in our result, such an invariant measure is indeed
absolutely continuous with respect to the Lebesgue measure.

This paper is organized as follows: in Section 2,  after introducing the notion of
almost everywhere stochastic (invertible) flow, we give two direct consequences of this notion and
then state our main results. In Section 3, we give some necessary preliminaries for later use. In Section 4,
we study stochastic transport equation (\ref{STR}) in case that
$b\in \mathrm{BV}_{loc}$ has bounded divergence
and $\sigma$ is constant. In Section 5, we apply the results of Section 4 to the study of
stochastic flows of SDE with BV drift and constant diffusion coefficients.
In Section 6, we extend the result of \cite{Cr-De-Le} to the stochastic case. Here, an SDE with discontinuous
coefficients is provided to show our result.
This section can be read independently of Sections 4 and 5. In Section 7, we prove our main results.
In the appendix, we give a detailed proof about the flow property as well as the Markov property
when SDE (\ref{SSDE}) admits a unique almost everywhere stochastic flow in the sense of Definition \ref{Def1} below.

\section{Main Results}

We first introduce some necessary notations. Let $(E,\cE, \mu)$ be a measure space and
$\sT:E\to E$ a measurable transformation. We shall use $\mu\circ\sT$ to denote the image
measure of $\mu$ under $\sT$, i.e., for any nonnegative measurable function $\varphi$,
$$
\int_E\varphi(x)\mu\circ\sT(\dif x):=\int_E\varphi(\sT(x))\mu(\dif x).
$$
By $\mu\circ\sT\ll\mu$ we mean that $\mu\circ\sT$ is absolutely continuous with respect to
$\mu$. Let $C^\infty_c(\mR^d)$ be the set of all smooth functions on $\mR^d$  with compact supports,
$C_b(\mR^d)$ the set of all bounded continuous functions,
and $\cL^+(\mR^d)$ the set of all nonnegative Borel measurable functions.
Below, we shall denote the Lebesgue measure by $\sL(\dif x)$ or $\dif x$.

\textsc{Convention:}
The repeated indices will be summed. The letter $C$ with or without subscripts
will denote a positive constant whose value is not important and may change in different occasions.
Moreover, all the derivatives, gradients and divergences are taken in the distributional sense.

We  introduce the following notion of almost everywhere stochastic (invertible) flows,
which is inspired by LeBris and Lions \cite{Le-Li} and Ambrosio \cite{Am}.
\bd\label{Def1}
Let $X_t(\omega,x)$ be a $\mR^d$-valued measurable
stochastic field on $\mR_+\times\Omega\times\mR^d$.
We say $X$ an {\bf almost everywhere stochastic flow} of (\ref{SSDE}) corresponding to $(b,\sigma)$ if
\begin{enumerate}[{\bf(A)}]
\item For  $\sL$-almost all $x\in\mR^d$, $t\mapsto X_t(x)$ is a continuous
($\cF_t$)-adapted stochastic process satisfying that for any $T>0$
$$
\int^T_0|b(X_s(x))|\dif s+\int^T_0|\sigma(X_s(x))|^2\dif s<+\infty, \ \ P-a.s.,
$$
and solves
$$
X_t(x)=x+\int^t_0b(X_s(x))\dif s+\int^t_0\sigma(X_s(x))
\dif W_s,\ \ \forall t\geq 0.
$$
\item For any $t\geq 0$ and $P$-almost all $\omega\in\Omega$,
$\sL\circ X_t(\omega,\cdot)\ll\sL$. Moreover, for any $T>0$, there exists a
constant $K_{T,b,\sigma}>0$ such that for all $\varphi\in\cL^+(\mR^d)$
\be\label{Den}
\sup_{t\in[0,T]}\mE\int_{\mR^d}\varphi(X_t(x))\dif x\leq K_{T,b,\sigma}
\int_{\mR^d}\varphi(x)\dif x.
\ee
\end{enumerate}
We say $X$ an {\bf almost everywhere stochastic  invertible flow} of (\ref{SSDE})
corresponding to $(b,\sigma)$ if
in addition to the above {\bf (A)} and {\bf (B)},

\begin{enumerate}[{\bf (C)}]
\item  For any $t\geq 0$ and $P$-almost all $\omega\in\Omega$, there exists a measurable
inverse $X^{-1}_t(\omega,\cdot)$ of $X_t(\omega,\cdot)$ so that
$\sL\circ X^{-1}_t(\omega,\cdot)=\rho_t(\omega,\cdot)\sL$, where the density $\rho_t(x)$ is given by
\be
\rho_t(x):=\exp\left\{\int^t_0\Big[\div b
-\frac{1}{2}\p_i\sigma^{jl}\p_j\sigma^{il}\Big](X_s(x))\dif s+\int^t_0
\div \sigma(X_s(x))\dif W_s\right\}.\label{Rho}
\ee
Here, $\div\sigma^{\cdot l}:=\p_i\sigma^{il}$ and we require that for
any $T>0$ and $\sL$-almost all $x\in\mR^d$,
$$
\int^T_0\Big[|\div b|+|\p_i\sigma^{jl}\p_j\sigma^{il}|+|\div\sigma|^2\Big]
(X_s(x))\dif s<+\infty,\ \ P-a.s.
$$
\end{enumerate}
\ed
\br
If $\sigma=constant$ and $\div b\in L^\infty(\mR^d)$, then {\bf (C)} clearly implies {\bf (B)}.
In fact, in this case we have
$$
\sL\circ X_t(\omega,\cdot)=\rho^{-1}_t(\omega,X^{-1}_t(\omega,\cdot))\sL
$$
and by (\ref{Rho})
$$
|\rho^{-1}_t(\omega,X^{-1}_t(\omega,x))|\leq e^{t\|\div b\|_\infty}.
$$
\er

In what follows, for the simplicity of notations,
we shall drop the time variable $t$ and the spatial variable $x$ if there are no confusions.
For examples, for a function $f_s(x)$, we simply write
$$
\iint f:=\iint_{\mR^d} f_s(x)\dif x\dif s
$$
and
$$
\iint f\dif W_s:=\iint_{\mR^d} f_s(x)\dif x\dif W_s.
$$

The following result is an easy consequence of Definition \ref{Def1}.
\bp\label{Pro1}
Assume that $b\in L^1_{loc}(\mR^d)$ with $\div b\in L^1_{loc}(\mR^d)$ and $\sigma\in C^2(\mR^d)$.
Let $X$ be an almost everywhere stochastic invertible flow of (\ref{SSDE}) in the sense of Definition \ref{Def1}.
Let $u_0\in L^\infty(\mR^d)$ and set $u_{t}(x):=u_0(X^{-1}_{t}(x))$.
Then $u_t(x)$ solves the following stochastic transport equation
in the distributional  sense:
$$
\dif u=\frac{1}{2}\sigma^{il}\sigma^{jl}\p^2_{ij}u
-b^i_\sigma\p_i u\dif t-\sigma^{il}\p_i u\dif W^l_t,
$$
where $b^i_\sigma:=b^i-\sigma^{jl}\p_j\sigma^{il}$.
In particular,  $\bar u_t(x):=\mE u_0(X^{-1}_t(x))$
is a distributional  solution of the following second order parabolic differential equation:
$$
\p_t \bar u=\frac{1}{2}\sigma^{il}\sigma^{jl}\p^2_{ij}\bar u
-b^i_\sigma\p_i \bar u.
$$
\ep
\begin{proof}
Let $\varphi\in C^\infty_c(\mR^d)$. By {\bf (C)} of Definition \ref{Def1}, we have
\ce
-\iint(b^i_\sigma\p_iu)\varphi&=&
\iint u_0(X^{-1})\cdot\div(b_\sigma\varphi)
=\iint u_0\cdot\div(b_\sigma\varphi)(X)\cdot\rho=\\
&=&\iint u_0\cdot(b^i_\sigma\p_i\varphi)(X)\cdot\rho
+\iint u_0\cdot (\varphi\div b_\sigma)(X)\cdot\rho.
\de
Similarly,
\ce
-\iint(\sigma^{il}\p_i u)\varphi \dif W^l_s
=\iint u_0\cdot(\sigma^{il}\p_i\varphi)(X)\cdot\rho \dif W^l_s
+\iint u_0\cdot(\p_i\sigma^{il}\varphi)(X)\cdot\rho \dif W^l_s
\de
and
\ce
\frac{1}{2}\iint\sigma^{il}\sigma^{jl}\p^2_{ij}u \varphi=
\frac{1}{2}\iint u_0\cdot[\p^2_{ij}(\sigma^{il}\sigma^{jl}\varphi)](X)\cdot\rho.
\de
Moreover, by stochastic Fubini's theorem, we have
\ce
&&\iint u_0\cdot (b^i_\sigma\p_i\varphi)(X)\cdot\rho
+\iint u_0\cdot (\sigma^{il}\p_i\varphi)(X)\cdot\rho\dif W^l_s\\
&&\quad=\int u_0\left(\int^t_0 (b^i_\sigma\p_i\varphi)(X)\cdot\rho\dif s
+\int^t_0(\sigma^{il}\p_i\varphi)(X)\cdot\rho\dif W^l_s\right)
\de
and
\ce
&&\iint u_0\cdot(\div b_\sigma\cdot\varphi)(X)\cdot\rho
+\iint u_0\cdot(\div \sigma\cdot\varphi)(X)\cdot\rho \dif W_s\\
&&\quad=\int u_0\left(\int^t_0(\div b_\sigma\cdot\varphi)(X)\cdot\rho\dif s
+\int^t_0(\div \sigma\cdot\varphi)(X)\cdot\rho\dif W_s\right).
\de

On the other hand, by (\ref{Rho}) and It\^o's formula, we have
\ce
\rho_t=1+\int^t_0\rho_s\Big[\div b_\sigma+\frac{1}{2}\p^2_{ij}
(\sigma^{il}\sigma^{jl})\Big](X_s)\dif s
+\int^t_0\rho_s \p_i\sigma^{i l}(X_s)\dif W^l_s,
\de
and
\ce
\dif[\varphi(X_t)\rho_t]&=&\Big[b^i\p_i\varphi
+\frac{1}{2}\sigma^{ik}\sigma^{jk}\p^2_{ij}\varphi\Big](X_t)\rho_t\dif t
+(\sigma^{il}\p_i\varphi)(X_t)\rho_t\dif W^l_t\\
&&+\Big[\varphi\div b_\sigma+\frac{1}{2}\varphi\p^2_{ij}
(\sigma^{il}\sigma^{jl})\Big](X_t)\rho_t\dif t
+[\varphi\p_i\sigma^{il}](X_t)\rho_t\dif W^l_t\\
&&+[\sigma^{il}\p_i\varphi\p_j\sigma^{jl}](X_t)\rho_t\dif t\\
&=&b^i_\sigma\p_i\varphi(X_t)\rho_t\dif t
+(\sigma^{il}\p_i\varphi)(X_t)\rho_t\dif W^l_t\\
&&+\varphi\div b_\sigma(X_t)\rho_t\dif t
+[\varphi\p_i\sigma^{il}](X_t)\rho_t\dif W^l_t\\
&&+\frac{1}{2}\p_{ij}^2(\sigma^{il}\sigma^{jl}\varphi)(X_t)\rho_t\dif t.
\de
Combining the above calculations, we get
\ce
&&\frac{1}{2}\iint\sigma^{il}\sigma^{jl}\p^2_{ij}u\varphi
-\iint(b^i_\sigma\p_i u)\varphi -
\iint(\sigma^{il}\p_i u)\varphi \dif W^l_s\\
&&\qquad=\int u_0 \left(\int^t_0\dif[\varphi(X_s)\rho_s]\right)
=\int u_0 [\varphi(X_t )\rho_t -\varphi ]\\
&&\qquad=\int u_0(X^{-1}_t)\varphi-\int u_0\varphi
=\int u_t\varphi-\int u_0\varphi.
\de
The proof is complete.
\end{proof}
The following proposition is much technical. We shall prove it in the appendix.
\bp\label{Pr2}
Assume that SDE (\ref{SSDE}) admits a unique almost everywhere stochastic (or invertible) flow.
Then the following flow property holds: for any $s\geq 0$ and $(P\times\sL)$-almost all
$(\omega,x)\in\Omega\times\mR^d$,
\be
X_{t+s}(\omega,x)=X_t(\theta_s\omega,X_s(\omega,x)),\  \forall t\geq 0,\label{Ep3}
\ee
where $\theta_s\omega:=\omega(s+\cdot)-\omega(s)$. Moreover, for any bounded measurable
function $\varphi$ on $\mR^d$, define
$$
\mT_t\varphi(x):=\mE\varphi(X_t(x)),
$$
then for any $t,s\geq 0$
\be
\mE(\varphi(X_{t+s}(x))|\cF_s)=\mT_t\varphi(X_s(x)),\ \ (P\times\sL)-a.e.\label{Ep30}
\ee
In particular, $(\mT_t)_{t\geq 0}$ forms a bounded linear operator semigroup
on $L^p(\mR^d)$ for any $p\geq 1$.
\ep
\br
Here, an open question is that whether the following stronger flow property holds:
For $(P\times\sL)$-almost all
$(\omega,x)\in\Omega\times\mR^d$,
\be
X_{t+s}(\omega,x)=X_t(\theta_s\omega,X_s(\omega,x)),\
\forall t,s\geq 0.\label{Ep33}
\ee
In the language of random dynamical systems (cf. \cite[Definition 1.1.1]{Ar}),
property (\ref{Ep3}) is called ``crude'', and property (\ref{Ep33}) is called ``perfect''.
A deep result of Arnold and Scheotzow (cf. \cite[p.17, Theorem 1.3.2]{Ar})
asserted that a crude cocycle admits an indistinguishable and perfect version.
But, it seems that we can not use
their result to deduce (\ref{Ep33}) since it is not clear how to endow a structure on
the set of all measurable transformations so that it becomes a Hausdorff topological group
with countable topological base.
\er

Our  main result of this paper is:
\bt\label{Main2}
Assume that
\be
\frac{|b(x)|}{1+|x|},\ \div b(x)\in L^\infty(\mR^d)\label{BB}
\ee
and one of the following conditions holds:
\be
&&b(x)\in \mathrm{BV}_{loc}\mbox{ and }\sigma\mbox{ is independent of $x$};\label{BB0}\\
&&\left\{ \begin{aligned}
&|\nabla b(x)|\in (L\log L)_{loc}(\mR^d), \\
&|\nabla\sigma(x)|,\ \ \sup_{|z|\leq 1}|\sigma(x-z)|\cdot|\nabla\div\sigma|(x)\in L^\infty(\mR^d).
\end{aligned} \right.\label{SI}
\ee
Then there exists a unique almost everywhere stochastic invertible flow of (\ref{SSDE}) corresponding to
$(b,\sigma)$ in the sense of Definition \ref{Def1}.
\et
\br
By definitions, $b\in BV_{loc}$ means that $\nabla b$ is a
locally finite vector valued Radon measure on $\mR^d$; and
$|\nabla b|\in (L\log L)_{loc}(\mR^d)$ means that $|\nabla b|\log(|\nabla b|+1)\in L^1_{loc}(\mR^d)$.
In particular, for any $p>1$,
$$
L^p_{loc}(\mR^d)\subset (L\log L)_{loc}(\mR^d)\subset L^1_{loc}(\mR^d).
$$
In (\ref{SI}), the second condition on $\sigma$ is certain growth restriction
of $\sigma$ and $\nabla\div\sigma$.
\er

About the existence of invariant measure of $(\mT_t)_{t\geq 0}$, we have the following criterion.
\bt\label{Main1}
Assume that SDE (\ref{SSDE}) admits a unique almost everywhere stochastic flow
with $K_{T,b,\sigma}=K_{b,\sigma}$ in (\ref{Den}) independent of $T$, and $(b,\sigma)$ satisfies
\be
\<x,b(x)\>_{\mR^d}+\|\sigma(x)\|_{H.S.}^2\leq 0( \mbox{ or $-C_1|x|^2+C_2$}),\label{CO}
\ee
where $C_1,C_2>0$, and $\|\sigma(x)\|_{H.S.}$ denotes the Hilbert-Schmidt norm of matrix $\sigma(x)$.
Then $(\mT_t)_{t\geq 0}$ admits an invariant probability measure $\mu(\dif x)=\gamma(x)\dif x$
with $\gamma\in L^\infty(\mR^d)\cap L^1(\mR^d)$ so that
for all $\varphi\in L^1(\mR^d)$ and $t\geq 0$
\be
\int_{\mR^d} \mT_t\varphi(x)\gamma(x)\dif x=\int_{\mR^d}\varphi(x)\gamma(x)\dif x.\label{Eq}
\ee
\et
\br
It is well known that if $\mT_t$ is a Feller semigroup, then under (\ref{CO}), there exists an
invariant probability measure for $\mT_t$. In our case, $\mT_t$ may be not a Feller semigroup.
In Theorem \ref{Th11} below, we shall give a condition such that $K_{T,b,\sigma}=K_{b,\sigma}$
in (\ref{Den}) is independent of $T$.
\er

These two theorems will be proved in Section 7.

\section{Preliminaries}

In this section, we prepare some lemmas for later use.
Below, we consider SDE (\ref{SSDE}) and assume that $b, \sigma\in C^\infty_b(\mR^d)$
are $C^\infty$-smooth, which together with their derivatives of all orders are bounded.
It is well known that
the family of solutions $\{X_t(x), t\geq 0\}_{x\in\mR^d}$ to SDE (\ref{SSDE}) forms a
$C^\infty$-diffeomorphism flow (cf. \cite{Ik-Wa, Ku}). We have the following simple result about
the Jacobian determinant of stochastic flow.
\bl\label{Le5}
Let $\rho_t(x)$ be defined by (\ref{Rho}). Then
\be
\det(\nabla X_t(x))=\rho_t(x)\label{Es4}
\ee
and for any $T>0$ and $p\geq 1$,
\be
\mE|\det(\nabla X^{-1}_T(x))|^p\leq\exp\left\{pT\Big(\|[-\div b+\frac{1}{2}\p_i\sigma^{jl}\p_j\sigma^{il}
+\sigma^{il}\p_{ij}^2\sigma^{jl}+\frac{p}{2}|\div \sigma|^2]^+\|_\infty\Big)\right\},\label{Es5}
\ee
where for a real number $a$, $a^+:=a\vee 0:=\max(a,0)$.
\el
\begin{proof}
Let $\tilde b^i:=b^i-\frac{1}{2}\sigma^{jl}\p_j\sigma^{il}$.
We write equation (\ref{SSDE}) as Stratonovich form:
$$
\dif X=\tilde b(X)\dif t+\sigma(X)\circ\dif W_t,\ \ X_0=x.
$$
Let $W^n_{t}$ be the linearized approximation of $W_t$. Consider the following ODE:
$$
\dif X_n(x)=\tilde b(X_n)\dif t+\sigma(X_n) \dot W^n_t\dif t.
$$
Then,
$$
\det(\nabla X_{n,t}(x))=\exp\left\{\int^t_0\div \tilde b(X_{n,s}(x))\dif s+\int^t_0
\div \sigma(X_{n,s}(x))\dot W^n_s\dif s\right\}.
$$
By the limit theorem (cf. \cite{Ik-Wa,Ku}), we get
$$
\det(\nabla X_t(x))=\exp\left\{\int^t_0\div \tilde b(X_s(x))\dif s+\int^t_0
\div \sigma(X_s(x))\circ \dif W_s\right\}.
$$
(\ref{Es4}) then follows by rewriting the Stratonovich integral as It\^o's integral.

On the other hand, fix $T>0$ and let $Y_s$ solve the following SDE:
$$
\dif Y_t=-\tilde b(Y_t)\dif t+\sigma(Y_t)\circ\dif W^T_t,\ \ Y_0=x,
$$
where $W^T_t:=W_{T-t}-W_T$. It is well known that (cf. \cite{Ik-Wa,Ku})
$$
X_T^{-1}(x)=Y_T(x).
$$
As above, we have
\ce
\det(\nabla Y_T)&=&\exp\left\{-\int^T_0\div \tilde b(Y_s)\dif s+\int^T_0
\div \sigma(Y_s)\circ \dif W^T_s\right\}\\
&=&\exp\left\{\int^T_0\Big[-\div \tilde b+\frac{1}{2}
\sigma^{il}\p^2_{ij}\sigma^{jl}\Big](Y_s)\dif s+\int^T_0
\div \sigma(Y_s)\dif W^T_s\right\}.
\de
Note that for any $p\geq 1$
$$
t\mapsto\exp\left\{p\int^t_0\div \sigma(Y_s)\dif W^T_s
-\frac{p^2}{2}\int^t_0|\div \sigma(Y_s)|^2\dif s\right\}
$$
is a continuous exponential martingale. Estimate (\ref{Es5}) then follows by H\"older's inequality.
\end{proof}

Let $C^\infty_p(\mR^d)$ be the set of all smooth functions with polynomial growth.
The following proposition is an easy consequence of Proposition \ref{Pro1}
(see also \cite[p.180, Theorem 1]{Ro}).
\bp\label{Pro25}
For any  $u_0\in C^\infty_p(\mR^d)$, let $u_t(x):=u_0(X^{-1}_t(x))$.
Then $u_t(x)$ solves the following
stochastic transport equation in the classical sense:
$$
\dif u=\frac{1}{2}\sigma^{il}\sigma^{jl}\p^2_{ij}u\dif t
-b^i_\sigma\p_i u\dif t-\sigma^{il}\p_i u\dif W^l_t,\ \ u|_{t=0}=u_0,
$$
where $b^i_\sigma:=b^i-\sigma^{jl}\p_j\sigma^{il}$.
\ep
The following result can be found in \cite{Ku} and \cite[p. 180, Theorem 1]{Ro}.
\bp\label{Pro26}
Let $X_{s,t}(x)$ solve
$$
X_{s,t}(x)=x+\int^t_sb(X_{s,r})\dif r+\int^t_s\sigma(X_{s,r})\dif W_r,\ \ t\geq s\geq0.
$$
Fix $t>0$. For any $v_0\in C^\infty_p(\mR^d)$, let
$v_{s,t}(x):=v_0(X_{s,t}(x))$, where $s\in[0,t]$.
Then $v_{s,t}(x)$ solves the following backward
stochastic Kolmogorov equation in the classical sense:
$$
\dif v+\frac{1}{2}\sigma^{il}\sigma^{jl}\p^2_{ij}v\dif s+b^i\p_i v\dif s
+\sigma^{il} \p_iv*\dif W_s=0,\ \ v|_{s=t}=v_0,
$$
where the asterisk denotes the backward It\^o's integral.
\ep

Let $C^+_c(\mR^d)$ be the set of all non-negative continuous functions on $\mR^d$
with compact support and $\sC$ a countable and dense subset of $C^+_c(\mR^d)$
with respect to the uniform norm $\|\varphi\|_\infty:=\sup_{x\in\mR^d}|\varphi(x)|$.
We need the following simple lemma.
\bl\label{Le3}
Let $X, Y:\mR^d\to\mR^d$ be two measurable transformations.

(i) Let $\gamma\in\cL^+(\mR^d)\cap L^1_{loc}(\mR^d)$. Assume that for any $\varphi\in\sC$,
\be
\int \varphi(X)\leq \int \varphi\cdot\gamma.\label{L22}
\ee
Then this inequality still holds for all $\varphi\in\cL^+(\mR^d)$. In particular, $\sL\circ X\ll\sL$.

(ii) Let $\rho:\mR^d\to\mR^+$ be a positive measurable function with $\rho\in L^1_{loc}(\mR^d)$.
Assume that for any $\varphi,\psi\in \sC$,
\be
\int \varphi(Y)\cdot\psi=\int \varphi\cdot\psi(X)\cdot\rho.\label{L11}
\ee
Then $X$ admits a measurable invertible $Y$, i.e., $X^{-1}(x)=Y(x)$ a.e.. Moreover,
$$
\sL\circ X^{-1}=\rho\sL,\ \ \ \sL\circ X=\rho^{-1}(X^{-1})\sL.
$$
\el
\begin{proof}
(i) Thanks to the density of $\sC$ in $C^+_c(\mR^d)$, by Fatou's lemma and the dominated convergence
theorem, one sees that (\ref{L22}) holds for all $\varphi\in C^+_c(\mR^d)$. Now, let $O\subset\mR^d$
be a bounded open set. Define
$$
\varphi_n(x):=1-\left(\frac{1}{1+\mathrm{distance}(x,O^c)}\right)^n.
$$
Then $\varphi_n\in C^+_c(\mR^d)$ and for every $x\in\mR^d$,
$$
\varphi_n(x)\uparrow 1_O(x)\mbox{ as $n\to\infty$}.
$$
By the monotone convergence theorem, we find that (\ref{L22}) holds for $\varphi=1_O$.
Thus, the desired conclusion follows by the monotone class theorem.

(ii)  As above, one sees that (\ref{L11})
holds for all $\varphi,\psi\in\cL^+(\mR^d)$. Thus, we have for all $\varphi,\psi\in\cL^+(\mR^d)$,
$$
\int \varphi(X\circ Y)\cdot\psi=\int \varphi(X)\cdot\psi(X)\cdot\rho=\int \varphi\cdot\psi
$$
and
$$
\int \varphi\cdot\psi(Y\circ X)\cdot\rho=\int \varphi(Y)\cdot\psi(Y)=\int \varphi\cdot\psi\cdot\rho.
$$
By the monotone class theorem, we obtain that for any Borel measurable
set $A\subset \mR^d\times\mR^d$,
$$
\int_{\mR^d} 1_A(X\circ Y(x),x)\cdot e^{-|x|}\dif x=\int_{\mR^d} 1_A(x,x)\cdot e^{-|x|} \dif x
$$
and
$$
\int_{\mR^d} 1_A(x,Y\circ X(x))\cdot e^{-|x|}\dif x=\int_{\mR^d} 1_A(x,x)\cdot e^{-|x|}\dif x.
$$
Hence, letting $A=\{(x,y):x\not=y\}$ yields that
$X\circ Y(x)=x$ and $Y\circ X(x)=x$ for $\sL$-almost all $x\in\mR^d$. The result follows.
\end{proof}

The following lemma will play a crucial role for taking limits below.
\bl\label{Le4}
Let $X_n(\omega,x):\Omega\times\mR^d\to\mR^d, n\in\mN$
be a family of measurable mappings, which are uniformly bounded in $L^\infty_{loc}(\mR^d;
L^p(\Omega))$ for any $p\geq 1$. Suppose that for $P$-almost all $\omega\in\Omega$,
$\sL\circ X_n(\omega,\cdot)\ll\sL$
and the density $\gamma_n(\omega,x)$ satisfies
\be
\sup_n~\mathrm{ess.}\sup_{x\in\mR^d}\mE|\gamma_n(x)|^2\leq C_1.\label{Es66}
\ee
If for ($P\times\sL$)-almost all $(\omega,x)\in\Omega\times\mR^d$,
$X_n(\omega,x)\to X(\omega,x)$ as
$n\to\infty$, then for $P$-almost all $\omega\in\Omega$, $\sL\circ X(\omega,\cdot)\ll \sL$ and
the density $\gamma$ also satisfies
\be
\mathrm{ess.}\sup_{x\in\mR^d}\mE|\gamma(x)|^2\leq C_1.\label{Es77}
\ee
Moreover, let $(\psi_n)_{n\in\mN}$ be a family of measurable functions
on $\mR^d$ and satisfy that for some
$C_2>0$ and $\a\geq 1$
\be
\sup_{n\in\mN}~\mathrm{ess.}\sup_{x\in\mR^d}\frac{|\psi_n(x)|}{1+|x|^\a}\leq C_2.\label{Ps1}
\ee
If $\psi_n$ converges to some $\psi$ in $L^1_{loc}(\mR^d)$, then for any $N>0$,
\be
\lim_{n\to\infty}\mE\int_{B_N}|\psi_n(X_n)-\psi(X)|=0.\label{Ps2}
\ee
\el
\begin{proof}
Fix $\varphi\in \sC\subset C^+_c(\mR^d)$ with support contained in $B_N$ for some $N>0$.
Then by Fubini's theorem and Fatou's lemma, we have for $P$-almost all $\omega\in\Omega$,
\be
\int \varphi(X(\omega))\leq\varliminf_{n\to\infty}\int \varphi(X_n(\omega))
=\varliminf_{n\to\infty}\int \varphi\cdot\gamma_{n}(\omega)=:
\varliminf_{n\to\infty}J^\varphi_n(\omega).
\label{PP2}
\ee
By (\ref{Es66}), there exists a subsequence still denoted by $n$ and
a $\gamma_0\in L^\infty(\mR^d; L^2(\Omega))$ satisfying (\ref{Es77}) such that
$$
\gamma_n\mbox{ weakly * converges to $\gamma_0$ in $L^\infty(\mR^d; L^2(\Omega))$}.
$$
Since $\gamma_n$ also weakly converges to $\gamma_0$ in $L^2(B_N\times\Omega)$,
by Banach-Saks' theorem, there is  another subsequence still denoted by $n$ such that
its Ces\`aro mean $\bar \gamma_n:=\frac{1}{n}\sum_{k=1}^n \gamma_k$ strongly converges
to $\gamma_0$ in $L^2(B_N\times\Omega)$.
Thus, there is another subsequence still denoted by $n$
such that for $P$-almost all $\omega\in\Omega$,
$$
\bar\gamma_n(\omega)\stackrel{n\to\infty}{\longrightarrow}\gamma_0(\omega)\mbox{ in $L^2(B_N)$}.
$$
Hence,
$$
\bar J^\varphi_n(\omega):=\frac{1}{n}\sum_{k=1}^nJ^\varphi_k(\omega)
=\int \varphi\cdot\bar\gamma_{n}(\omega)
\stackrel{n\to\infty}{\longrightarrow}\int \varphi\cdot\gamma_0(\omega),
$$
which together with (\ref{PP2}) yields that for $P$-almost all $\omega$,
$$
\int \varphi(X(\omega))\leq \varliminf_{n\to\infty}J^\varphi_n(\omega)\leq
\lim_{n\to\infty}\bar J^\varphi_n(\omega)=\int \varphi\cdot\gamma_0(\omega).
$$
Since $\sC$ is countable, we may find a common null set $\Omega'\subset\Omega$ such that
the above inequality holds for all $\omega\notin\Omega'$ and $\varphi\in\sC$.
The first conclusion then follows by (i) of Lemma \ref{Le3}.

We now prove (\ref{Ps2}). We make the following decomposition:
\ce
&&\int_{B_N}|\psi_n(X_n)-\psi(X)|\leq\int_{B_N}|\psi_n(X_n)-\psi(X_n)|\\
&&\qquad+\int_{B_N}|\psi(X_n)-\psi(X)|=:I_n+J_n.
\de
By (\ref{Ps1}) and $\psi_n\to\psi$ in $L^1_{loc}(\mR^d)$, we also have
$$
\mathrm{ess.}\sup_{x\in\mR^d}\frac{|\psi(x)|}{1+|x|^\a}\leq C_2.
$$
Let $(\phi_m)_{m\in\mN}$ be a family of bounded continuous functions
such that $\phi_m\to\psi$ in $L^1_{loc}(\mR^d)$ as $m\to\infty$ and
\be
\sup_{m\in\mN}~\mathrm{ess.}\sup_{x\in\mR^d}\frac{|\phi_m(x)|}{1+|x|^\a}\leq C_2.\label{Ps11}
\ee
We have
\ce
J_n&\leq&\int_{B_N}|\phi_m(X_n)-\psi(X_n)|+\int_{B_N}|\phi_m(X)-\psi(X)|\\
&&+\int_{B_N}|\phi_m(X_n)-\phi_m(X)|=:J_{1nm}+J_{2m}+J_{3nm}.
\de
For any $R>0$, we may write
\ce
&&J_{1nm}=\int_{B_N\cap\{|X_n|\leq R\}}|\phi_m(X_n)-\psi(X_n)|\\
&&\qquad+\int_{B_N\cap\{|X_n|>R\}}|\phi_m(X_n)-\psi(X_n)|=:J_{1nm}^{1,R}+J^{2,R}_{1nm}.
\de
By the change of variable and (\ref{Es66}), we have
$$
\mE J_{1nm}^{1,R}\leq\int_{B_R}|\phi_m-\psi|\cdot\mE\gamma_n\leq C_1\int_{B_R}|\phi_m-\psi|.
$$
By Chebyshev's inequality and (\ref{Ps11}), we have
$$
\mE J^{2,R}_{2nm}\leq \frac{C_{N,\a}}{R}\sup_{x\in B_N}\mE(1+|X_n(x)|^{2\a})\leq \frac{C_{N,\a}}{R}.
$$
Combining the above two estimates, we obtain
\be
\lim_{m\to\infty}\sup_{n\in\mN}\mE J_{1nm}=0.\label{Lp6}
\ee
Similarly, we also have
$$
\lim_{m\to\infty}\mE J_{2m}=0
$$
and for fixed $m\in\mN$, by the dominated convergence theorem,
$$
\lim_{n\to\infty}\mE J_{3nm}=0.
$$
Hence,
$$
\lim_{n\to\infty}\mE J_{n}=0.
$$
As proving (\ref{Lp6}), we also have
$$
\lim_{n\to\infty}\mE I_{n}=0.
$$
The proof is then complete.
\end{proof}
The following lemma will be used to prove the strong convergence in Theorem \ref{Pro3} below.
\bl\label{Le7}
Let $\mB$ be a separable and uniformly convex Banach space. Let
$(u_n)_{n\in\mN}$ be a bounded sequence in $L^1(\Omega; C([0,T];\mB))$.
Assume that for some $u\in L^1(\Omega; C([0,T];\mB))$,
\be
\lim_{n\to\infty}\mE\left(\sup_{t\in[0,T]}|
{}_{\mB^*}\<\phi,u_n(t)-u(t)\>_\mB|\right)=0,\ \ \forall\phi\in\mB^*,\label{P4}
\ee
where $\mB^*$ is the dual space of $\mB$, and
\be
\lim_{n\to\infty}\mE\left(\sup_{t\in[0,T]}\big|\|u_n(t)\|_{\mB}-\|u(t)\|_\mB\big|\right)=0.\label{P44}
\ee
Then $\sup_{t\in[0,T]}\|u_n(t)-u(t)\|_{\mB}$ converges to zero in probability as $n\to\infty$.
\el
\begin{proof}
It is enough to prove that for any subsequence $n_k$, there exists a subsubsequence $n'_k$
such that $\sup_{t\in[0,T]}\|u_{n_k'}(t)-u(t)\|_{\mB}$ converges to zero $P$-almost surely as
$k\to\infty$. We now fix a subsequence $n_k$ below. Since $\mB^*$ is separable, by (\ref{P4})
and (\ref{P44}),
we may find a subsubsequence $n_k'$ and a  measurable set $\Omega'\subset\Omega$ with $P(\Omega')=1$
such that for all $\omega\in \Omega'$, $u(\omega,\cdot)\in C([0,T];\mB)$ and
\be
\lim_{k\to\infty}\sup_{t\in[0,T]}|{}_{\mB^*}\<\phi,u_{n'_k}(\omega,t)-u(\omega,t)\>_\mB|=0,\ \
\forall\phi\in\mB^*\label{LL3}
\ee
and
\be
\lim_{k\to\infty}\sup_{t\in[0,T]}\big|\|u_{n'_k}(\omega,t)\|_{\mB}-\|u(\omega,t)\|_\mB\big|=0.\label{LL4}
\ee
We want to show that for such $\omega\in\Omega'$,
$$
\lim_{k\to\infty}\sup_{t\in[0,T]}\|u_{n'_k}(\omega,t)-u(\omega,t)\|_\mB=0.
$$
Suppose that this is not true. Then, there exist a $\delta>0$
and a sequence $(t_k)_{k\in\mN}\subset[0,T]$ such that
\be
\|u_{n'_k}(\omega,t_k)-u(\omega,t_k)\|_\mB\geq\delta,\ \ \forall k\in\mN.\label{LL5}
\ee
Without loss of generality, we assume that $t_k$ converges to $t_0$.
By (\ref{LL3}), (\ref{LL4}) and $u(\omega,\cdot)\in C([0,T];\mB)$,
we have
$$
\lim_{k\to\infty}\|u_{n'_k}(\omega,t_k)-u(\omega,t_0)\|_\mB=0,
$$
which together with $u(\omega,\cdot)\in C([0,T];\mB)$ yields
$$
\lim_{k\to\infty}\|u_{n'_k}(\omega,t_k)-u(\omega,t_k)\|_\mB=0.
$$
This is a contradiction with (\ref{LL5}). The proof is complete.
\end{proof}
We also recall some facts about local maximal functions. Let $f$ be a locally integrable function
on $\mR^d$. For every $R>0$, the local maximal function is defined by
$$
M_R f(x):=\sup_{0<r<R}\frac{1}{|B_r|}\int_{B_r} f(x+y)\dif y
=:\sup_{0<r<R}\fint_{B_r} f(x+y)\dif y.
$$
The following result can be found in \cite[p.143, Theorem 3]{Ev-Ga} and \cite[Appendix A]{Cr-De-Le}.
\bl

(i) (Morrey's inequality)
Let $f\in L^1_{loc}(\mR^d)$ be such that $\nabla f\in L^q_{loc}(\mR^d)$ for some $q>d$. Then
there exist $C_{q,d}>0$ and a negligible set $A$ such that
for all $x,y\in A^c$ with $|x-y|\leq R$,
\be
|f(x)-f(y)|&\leq& C_{q,d}\cdot |x-y|\cdot\left(\fint_{B_{|x-y|}}
|\nabla f|^q(x+z)\dif z\right)^{1/q}\no\\
&\leq& C_{q,d}\cdot |x-y|\cdot (M_R|\nabla f|^q(x))^{1/q}.\label{Es02}
\ee
(ii) Let $f\in L^1_{loc}(\mR^d)$ be such that $\nabla f\in L^1_{loc}(\mR^d)$.
Then there exist $C_d>0$ and a negligible set $A$ such that
for all $x,y\in A^c$ with $|x-y|\leq R$,
\be
|f(x)-f(y)|\leq C_d\cdot |x-y|\cdot(M_R |\nabla f|(x)+M_R|\nabla f|(y)).\label{Es2}
\ee
(iii) Let $f\in (L\log L)_{loc}(\mR^d)$. Then for any $N,R>0$ and some $C_{d,N}, C_{d}>0$,
\be
\int_{B_N}M_R|f|\leq C_{d,N}+C_{d}\int_{B_{N+R}}|f|\log(|f|+1).\label{Es3}
\ee
(iv) Let $f\in L^p_{loc}(\mR^d)$ for some $p>1$. Then for some $C_{d,p}>0$ and any $N,R>0$,
\be
\left(\int_{B_N}(M_R|f|)^p\right)^{1/p}\leq C_{d,p}\left(\int_{B_{N+R}}|f|^p\right)^{1/p}.\label{Es30}
\ee
\el

\section{Stochastic Transport Equations}

In this section we work on $[0,T]$ and mainly study the following stochastic transport equation:
\be
\dif u=\Big[\frac{1}{2}\sigma^{il}\sigma^{jl}\p^2_{ij}u+b^i\p_iu\Big]\dif t
+\sigma^{il}\p_iu\dif W^l_t,\ \ u|_{t=0}=u_0,\label{Eq11}
\ee
where $\sigma\in\mR^d\times\mR^m$ does  not depend on $x$, and $b$ is a BV vector field and satisfies
\be
\frac{b(x)}{1+|x|},\ \div b(x)\in L^\infty(\mR^d),\ \ b\in \mathrm{BV}_{loc}.\label{Con}
\ee

We first introduce the following notion of renormalized solutions for equation (\ref{Eq1}).
\bd\label{Def3}
A measurable and ($\cF_t$)-adapted stochastic field $u: [0,T]
\times\Omega\times\mR^d\to\mR$ is called a renormalized solution
of (\ref{Eq11}) if for any $\beta\in C^2(\mR)$,
$$
v_t(\omega,x):=\beta(\arctan u_t(\omega,x))
$$
solves (\ref{Eq11}) in the distributional  sense, i.e., for any $\phi\in C^\infty_c(\mR^d)$
\be
\int v_t\phi=\int v_0\phi+\frac{1}{2}\iint v\sigma^{il}
\sigma^{jl}\p^2_{ij}\phi-\iint v(\div b\phi+b^i\p_i\phi)
-\iint v\sigma^{il}\p_i\phi\dif W^l_s.\label{Gen}
\ee
\ed
\br
Since $v$ is bounded, it is clear that both sides of (\ref{Gen}) are well defined.
\er

Our main result in this section is that
\bt\label{Th1}
Assume that condition (\ref{Con}) holds.

(Existence and Uniqueness) For any measurable function $u_0$, there exists a unique renormalized solution $u$
to stochastic transport equation (\ref{Eq11}) with $u|_{t=0}=u_0$
in the sense of Definition \ref{Def3}. Moreover, for any $p>1$ and $N>0$,
$$
\arctan u\in L^p(\Omega; C([0,T];L^p(B_N))).
$$

(Stability) Let $b_n\in L^1_{loc}(\mR^d)$ be such that $\div b_n\in L^1_{loc}(\mR^d)$ and
$b_n, \div b_n$ converge to $b,\div b$ respectively in $L^1_{loc}(\mR^d)$.
Let $u^n_0$ $\sL$-almost everywhere converge to $u_0$. Let $u^n$ and $u$ be the renormalized
solutions corresponding to $(b^n, u^n_0)$ and $(b,u_0)$ in the sense of Definition \ref{Def3}.
Then for any $p>1$ and $N>0$,
$$
\arctan u^n\to \arctan u\ \ \mbox{strongly in $L^p(\Omega; C([0,T];L^p(B_N)))$}.
$$
\et

For proving this theorem, we first study the following more general stochastic
partial differential equation:
\be
\dif u=\Big[\frac{1}{2}\sigma^{il}\sigma^{jl}\p^2_{ij}u+b^i\p_iu+cu\Big]\dif t
+(\sigma^{il}\p_iu+h^lu)\dif W^l_t,\ \ u|_{t=0}=u_0,\label{Eq1}
\ee
where $\sigma$ and $b$ are as above and
\be
c, h, \nabla h\in L^\infty(\mR^d).\label{Con1}
\ee

As Definition \ref{Def3}, we also introduce the following notion about the renormalized solutions
for equation (\ref{Eq1}).
\bd\label{Def2}
We say $u\in L^\infty([0,T]\times\Omega\times\mR^d)$ a renormalized solution
of (\ref{Eq1}) if for any $\beta\in C^2(\mR)$, it holds that in the distributional  sense
\ce
\dif\beta(u)&=&\Big[\frac{1}{2}\sigma^{il}\sigma^{jl}\p^2_{ij}\beta(u)+b^i\p_i\beta(u)
+cu\beta'(u)\Big]\dif t\\
&&+\Big[\frac{1}{2}|h|^2\beta''(u)u^2+h^l\sigma^{il}\p_i\beta'(u)u\Big]\dif t\\
&&+(\sigma^{il}\p_i\beta(u)+h^lu\beta'(u))\dif W^l_t.
\de
\ed

We remark that for equation (\ref{Eq1}), the renormalized solution is a nonlinear notion,
whereas the distributional solution is a linear notion. However,
under (\ref{Con}) and (\ref{Con1}), we can show that these two notions are equivalent.
For this aim, we need the following class of regularized functions:
\be
\cN:=\left\{\varrho\in C^\infty_c(B_1),\ \ \varrho\geq 0,\ \ \int\varrho=1\right\}.\label{PL3}
\ee

We now establish the following equivalence between the distributional  solution and renormalized
solution.
\bp\label{Pro4}
Let $u\in L^\infty([0,T]\times\Omega\times\mR^d)$ be a distributional  solution of
(\ref{Eq1}). Then under (\ref{Con})  and (\ref{Con1}),
$u$ is also a renormalized solution of (\ref{Eq1}) in the sense of Definition \ref{Def2}.
\ep
\begin{proof}
Let $\varrho\in\cN$ and set $\varrho_\eps(x):=\eps^{-d}\varrho(x/\eps)$.
Define
$$
u_\eps:=u_{t,\eps}(x):=u_t*\varrho_\eps(x)=\int u_t(y)\varrho_\eps(x-y)\dif y.
$$
Taking convolutions for both sides of (\ref{Eq1}), we obtain
$$
\dif u_{\eps}= \Big[\frac{1}{2}\sigma^{il}\sigma^{jl}\p^2_{ij}u_{\eps}
+(b^i\p_iu)*\varrho_\eps+(cu)*\varrho_\eps\Big]\dif t
+[\sigma^{il}\p_iu_{\eps}+(h^lu)*\varrho_\eps]\dif W^l_t.
$$
Let $\beta\in C^2(\mR)$. By It\^o's formula, we have
\ce
\dif \beta(u_{\eps})&=&\Big[\frac{1}{2}\sigma^{il}\sigma^{jl}\p^2_{ij}\beta(u_\eps)
+((b^i\p_iu)*\varrho_\eps+(cu)*\varrho_\eps)\cdot\beta'(u_{\eps})\Big]\dif t\\
&&+\Big[\frac{1}{2}\beta''(u_{\eps})((h^lu)*\varrho_\eps)^2
+\beta''(u_{\eps})\sigma^{il}\p_iu_{\eps}\cdot(h^lu)*\varrho_\eps\Big]\dif t\\
&&+[\sigma^{il}\p_i\beta(u_{\eps})+(h^lu)*\varrho_\eps\cdot\beta'(u_{\eps})]\dif W^l_t.
\de
Write
$$
r^\rho_{\eps}:=((b^i\p_iu)*\varrho_\eps-b^i\p_i(u*\varrho_\eps))\cdot\beta'(u_{\eps})
$$
and
$$
[\varrho_\eps, h^l](u):=(h^lu)*\varrho_\eps-h^l(u*\varrho_\eps).
$$
Let $\phi\in C^\infty_c(\mR^d)$. Multiplying both sides by $\phi$ and integrating over $\mR^d$,
by the integration by parts formula, we get
\ce
\int\beta(u_{t,\eps})\phi&=&\int\beta(u_{0,\eps})\phi+
\iint\beta(u_\eps)\Big[\frac{1}{2}\sigma^{il}\sigma^{jl}\p^2_{ij}\phi
-\div b\phi-b^i\p_i\phi\Big]\\
&&+\iint r^\varrho_\eps\phi+\iint(cu)*\varrho_\eps)\cdot\beta'(u_{\eps})\phi\\
&&+\iint\frac{1}{2}\beta''(u_{\eps})((h^lu)*\varrho_\eps)^2\phi
+\iint\Big[\beta(u_\eps)-u_\eps\beta'(u_{\eps})\Big]\sigma^{il}\p_i(h^l\phi)\\
&&-\iint\beta'(u_\eps)\sigma^{il}\Big(\p_i[\varrho_\eps,h^l](u)\phi+[\varrho_\eps,h^l](u)\p_i\phi\Big)\\
&&+\iint[(h^lu)*\varrho_\eps\cdot\beta'(u_{\eps})\phi-\sigma^{il}\p_i\phi\beta(u_{\eps})]\dif W^l_s.
\de
Now taking limits $\eps\to 0$ and using \cite[p.516, Lemma II.1]{Di-Li}, we find that
\ce
&&\Bigg|\int\beta(u_{t})\phi-\int\beta(u_{0})\phi-
\iint\beta(u)\Big[\frac{1}{2}\sigma^{il}\sigma^{jl}\p^2_{ij}\phi
-\div b\phi-b^i\p_i\phi\Big]\\
&&-\iint(cu)\cdot\beta'(u)\phi-\iint\frac{1}{2}\beta''(u)(h^lu)^2\phi
-\iint\Big[\beta(u)-u_\eps\beta'(u)\Big]\sigma^{il}\p_i(h^l\phi)\\
&&-\iint[(h^lu)\cdot\beta'(u)\phi-\sigma^{il}\p_i\phi\beta(u)]\dif W^l_s\Bigg|\leq
\limsup_{\eps\to 0}\left|\iint r^\varrho_{\eps} \phi\right|.
\de
Since the left hand side of the above inequality does not depend on $\varrho$,
it suffices to show that
$$
\inf_{\varrho\in\cN}\limsup_{\eps\to 0}\left|\iint r^\varrho_{\eps} \phi\right|=0.
$$
This has been proved in the proof of \cite[Theorem 3.5]{Am}.
\end{proof}

Using Proposition \ref{Pro4}, we can prove the uniqueness of distributional  solutions.
\bp\label{Pro2}
Let $u\in L^\infty([0,T]\times\Omega; L^1(\mR^d)\cap L^\infty(\mR^d))$ be a distributional  solution
of (\ref{Eq1}). If $u|_{t=0}=0$, then
$$
u_t(\omega,x)=0,\ \ a.e.
$$
\ep
\begin{proof}
Let $\chi\in C^\infty_c(\mR^d)$ be a nonnegative cutoff function
with
\begin{equation} \label{eq:4}
\|\chi\|_\infty\leq 1,\ \
\chi(x)=\left\{ \begin{aligned}
&1,\ \ |x|\leq 1, \\
&0,\ \ |x|\geq 2.
\end{aligned} \right.
\end{equation}
Set $\chi_n(x):=\chi(x/n)$. By Proposition \ref{Pro4} and Definition \ref{Def2},
we have
\ce
\mE\int u^2_t\chi_n&=&\mE\iint\Big[\frac{1}{2}u^2\sigma^{il}\sigma^{jl}\p^2_{ij}\chi_n-u^2b^i\p_i\chi_n\Big]\\
&&+\mE\iint\Big[-u^2\div b\chi_n+2cu^2\chi_n\Big]\\
&&+\mE\iint\Big[|h|^2u^2\chi_n-u^2\p_i(h^l\chi_n)\sigma^{il}\Big].
\de
Observe that by (\ref{Con})
\be
|b^i\p_i\chi_n|\leq
\frac{|b|\cdot 1_{\{n\leq|x|\leq 2n\}}\cdot\|\nabla\chi\|_\infty}{n}\leq
C_1\cdot 1_{\{|x|\geq n\}},\label{EE2}
\ee
where $C_1=3\|b/(1+|x|)\|_\infty\cdot\|\nabla \chi\|_\infty$, and
$$
|h^l\p_i\chi_n|\leq \|h\|_\infty\cdot\|\nabla\chi\|_\infty.
$$
Since $u^2\in L^\infty([0,T]\times\Omega; L^1(\mR^d))$,
by letting $n\to\infty$, we obtain
\ce
\mE\int u_t^2&=&\mE\iint\Big[(-\div b+2c+|h|^2-\sigma^{il}\p_i h^l)u^2_s\Big]\\
&\leq&\|2c+|h|^2-\div b-\sigma^{il}\p_i h^l\|_\infty\int^t_0\left(\mE\int u^2_s\right)\dif s,
\de
which gives by Gronwall's inequality that
$$
\mE\int u_t^2=0.
$$
The uniqueness follows.
\end{proof}

In general, it is not expected to have a bounded solution for SPDE (\ref{Eq1}) because of the
presence of stochastic integral $\int^t_0h^l u\dif W^l_s$ (cf. \cite{Ro}).
We now turn back to stochastic transport equation (\ref{Eq11}), and
prove the existence-uniqueness and stability of $L^\infty$-distributional  solutions
when the initial value belongs to $L^\infty(\mR^d)$.
\bt\label{Pro3}
Assume that condition (\ref{Con}) holds.

(Existence and Uniqueness) For any $u_0\in L^\infty(\mR^d)$, there exists a
unique distributional  solution $u\in L^\infty([0,T]\times\Omega\times\mR^d)$
(also a renormalized solution in the sense of Definition \ref{Def2}) to
stochastic transport equation (\ref{Eq11}) satisfying
\be
\|u_t(\omega)\|_\infty\leq \|u_0\|_\infty.\label{Op1}
\ee
Moreover, there is a version still denoted by $u$ such that for any $p>1$ and $N>0$
\be
u\in L^p(\Omega; C([0,T];L^p(B_N))).\label{Lp5}
\ee

(Stability) Let $b_n\in L^1_{loc}(\mR^d)$ and $u^n_0\in L^\infty(\mR^d)$
be such that $\div b_n\in L^1_{loc}(\mR^d)$ and
$b_n, \div b_n, u^n_0$ converge to $b,\div b, u_0$ respectively in $L^1_{loc}(\mR^d)$.
Let $u_n, u\in L^\infty([0,T]\times\Omega\times\mR^d)$ be
the distributional  solutions of (\ref{Eq11}) corresponding to $(b_n, u^n_0)$ and $(b,u_0)$
and satisfy (\ref{Lp5}). Assume that
\be
\sup_n\|u_n\|_{L^\infty([0,T]\times\Omega\times\mR^d)}
<+\infty.\label{Ep1}
\ee
Then for any $p>1$ and $N>0$,
\be
u_n\to u\ \ \mbox{strongly in $L^p(\Omega; C([0,T]; L^p(B_N)))$}.\label{P8}
\ee
\et
\begin{proof}
(Existence) Fix a $\varrho\in\cN$ and a cutoff function $\chi$ satisfying
(\ref{eq:4}). Let
$$
\varrho_n(x):=n^d\varrho(nx),\ \ \chi_n(x)=\chi(x/n)
$$
and define
\be
b_n=b*\varrho_n\cdot\chi_n.\label{BN}
\ee
Let $X_n$ solve the following SDE:
$$
\dif X_n=-b_n(X_n)\dif t-\sigma\dif W_t,\  \ X_n|_{t=0}=x.
$$
By Proposition \ref{Pro25}, $u_{n,t}:=u_0(X^{-1}_{n,t})$ solves the following SPDE:
$$
\dif u_n=\Big[\frac{1}{2}\sigma^{il}\sigma^{jl}\p^2_{ij}u_n+b^i_n\p_iu_n\Big]\dif t
+\sigma^{il}\p_iu_n\dif W^l_t,\ \ u_n|_{t=0}=u_0.
$$
Clearly, $u_n\in L^\infty([0,T]\times\Omega\times\mR^d)$ and
$$
\|u_{n,t}(\omega,\cdot)\|_\infty\leq \|u_0\|_\infty.
$$
Therefore, for some $u\in L^\infty([0,T]\times\Omega\times\mR^d)$ and some subsequence
$n_k$,
$$
u_{n_k}\to u \mbox{ weakly* in $L^\infty([0,T]\times\Omega\times\mR^d)$}.
$$
Taking weakly* limits, it is easy to see that $u$ is a distributional  solution of (\ref{Eq11}).
Moreover, (\ref{Op1}) holds. As for (\ref{Lp5}), it can be seen from the
proof of the following stability.

(Uniqueness) Let $u$ and $\hat u$ be two distributional  solutions of (\ref{Eq11}) with the same
initial value. Then $v:=u-\hat u\in L^\infty([0,T]\times\Omega\times\mR^d)$ is
still a distributional  solution of (\ref{Eq11}) with zero
initial value. Since $v$ does not belong to
$L^\infty([0,T]\times\Omega;L^1(\mR^d))$, we can not directly use Proposition \ref{Pro2}
to obtain $v=0$. Below, we use a simple trick. Let
$$
\lambda(x):=\frac{1}{(1+|x|^2)^d},\ \  \hat v_t:=v_t\cdot\lambda.
$$
It is easy to see that
$$
\hat v_t\in L^\infty([0,T]\times\Omega;L^1(\mR^d)\cap L^\infty(\mR^d)).
$$
Moreover, noting that
$$
\p_i\lambda(x)=-\frac{2d x_i}{1+|x|^2}\lambda(x),\ \
\p_i\p_j\lambda(x)=\left(\frac{4d(d+1)x_ix_j}{(1+|x|^2)^2}-\frac{2d \delta_{ij}}{1+|x|^2}
\right)\lambda(x),
$$
we can check that $\hat v_t$ is  a distributional  solution of
$$
\dif \hat v=\Big[\frac{1}{2}\sigma^{il}\sigma^{jl}\p^2_{ij}\hat v+\hat b^i\p_i\hat v
+c\hat v\Big]\dif t
+(\sigma^{il}\p_i\hat v+h^l\hat v)\dif W^l_t,\ \ \hat v|_{t=0}=0,
$$
where
$$
\hat b^i(x)=b^i(x)+\frac{2dx_j\sigma^{il}\sigma^{jl}}{1+|x|^2},\ \
h^l(x)=\frac{2d x_i\sigma^{il}}{1+|x|^2}
$$
and
$$
c(x)=\frac{2d x_ib^i(x)}{1+|x|^2}
+\left(\frac{2d(d-1)x_ix_j}{(1+|x|^2)^2}+\frac{d \delta_{ij}}
{1+|x|^2}\right)\sigma^{il} \sigma^{jl}.
$$
By (\ref{Con}), one sees that $\hat b$ still satisfies (\ref{Con}) and
$c,h$ satisfy (\ref{Con1}).  Thus, we can use Proposition \ref{Pro2} to get
$\hat v=0$. The uniqueness follows.

(Stability) We follow DiPerna-Lions' argument \cite[p.523]{Di-Li}.
Fix an even number $p>1$ and let $v_n:=u_n^p$. Then,
by Definition \ref{Def2}, $v_n$ is  a distributional  solution of
$$
\dif v_n=\Big[\frac{1}{2}\sigma^{il}\sigma^{jl}\p^2_{ij}v_n+b^i_n\p_iv_n\Big]\dif t
+\sigma^{il}\p_iv_n\dif W^l_t,\ \ v_n|_{t=0}=(u^n_0)^p.
$$
By (\ref{Ep1}), we have
$$
\sup_n\|u_n\|_{L^\infty([0,T]\times\Omega\times\mR^d)}
+\sup_n\|v_n\|_{L^\infty([0,T]\times\Omega\times\mR^d)}
<+\infty.
$$
Without loss of generality, we may assume that $u_n$ and $v_n$ converges
weakly* in $L^\infty([0,T]\times\Omega\times\mR^d)$ to $u$ and $v$, which are
distributional  solutions of (\ref{Eq11}) corresponding to $u|_{t=0}=u_0$
and $v|_{t=0}=u_0^p$ by the assumptions.
By Proposition \ref{Pro4} and the uniqueness proved above, we have
$$
u^p=v.
$$
Thus,
$$
u_n^p\to u^p\ \ \mbox{weakly* in $L^\infty([0,T]\times\Omega\times\mR^d)$}.
$$
Hence, for any $N>0$,
$$
\mE\int^T_0\!\!\!\int_{B_N}u_n^p\to\mE\int^T_0\!\!\!\int_{B_N}u^p.
$$
By virtue of
$$
u_n\to u \ \ \mbox{weakly in $L^p([0,T]\times\Omega\times B_N)$},
$$
we thus obtain that for any $N>0$,
\be
u_n\to u\ \ \mbox{strongly in $L^p([0,T]\times\Omega\times B_N)$}.\label{P5}
\ee

We now strengthen this convergence to (\ref{P8}).
Let $w_n=u_n-u$. Then we have for any $N>0$ and $\phi\in C^\infty_c(B_N)$,
\ce
\left|\int w_{n,t}\phi\right|&=&\Bigg|\int w_{n,0}\phi+\iint\Big[w_n\frac{1}{2}\sigma^{il}\sigma^{jl}\p^2_{ij}\phi
-w_n(\div b\phi+b^i\p_i\phi)\Big]\\
&&-\iint u_n\Big[\div (b_n-b)\phi+(b^i_n-b^i)\p_i\phi\Big]
-\int^t_0\left(\int w_n\sigma^{il}\p_i \phi\right)\dif W^l_s\Bigg|\\
&\leq&C\int_{B_N}|w_{n,0}|+C\iint_{B_N}|w_n|+C\iint_{B_N}\Big[|\div (b_n-b)|+|b^i_n-b^i|\Big]\\
&&+\left|\int^t_0\left(\int w_n\sigma^{il}\p_i \phi\right)\dif W^l_s\right|.
\de
Hence, by BDG's inequality, (\ref{P5}) and the assumptions, we get
$$
\lim_{n\to\infty}\mE\left(\sup_{t\in[0,T]}\left|\int w_{n,t}\phi\right|\right)
\leq C\lim_{n\to\infty}\mE\left(\int^T_0\left(\int_{B_N} |w_n|\right)^2\dif s\right)^{1/2}=0.
$$
By another approximation, we further have for any $N>0$ and $\phi\in L^{p/(p-1)}(B_N)$,
\be
\lim_{n\to\infty}\mE\left(\sup_{t\in[0,T]}\left|\int_{B_N} w_{n,t}\phi\right|\right)=0.\label{P9}
\ee
Similarly, we also have for any $N>0$,
\be
\lim_{n\to\infty}\mE\left(\sup_{t\in[0,T]}\left|\int_{B_N}(u^p_{n,t}-u^p_t)\right|\right)=0.\label{P10}
\ee
Combining (\ref{Ep1}), (\ref{P9}) and (\ref{P10}), we then obtain (\ref{P8}) by Lemma \ref{Le7}.
\end{proof}

We are now in a position to give:

{\it Proof of Theorem \ref{Th1}:} (Uniqueness) Let $u$ and $\hat u$
be two renormalized solutions of SPDE (\ref{Eq11})
corresponding to the initial value $u_0$ in the sense of Definition \ref{Def3}.
Then $\arctan u$ and $\arctan\hat u$
are two distributional  solutions of SPDE (\ref{Eq11})
corresponding to the initial value $\arctan u_0$.
By Proposition \ref{Pro2}, we have
$$
\arctan u=\arctan \hat u.
$$
Hence,
$$
u=\hat u.
$$

(Existence) Let $v$ be the unique renormalized solution of SPDE (\ref{Eq11}) given in Proposition
\ref{Pro3} corresponding to the initial value $\arctan u_0\in[-\pi/2,\pi/2]$. Since
$$
\|v_t(\omega)\|_\infty\leq \|\arctan u_0\|_\infty\leq\pi/2,
$$
we may define
$$
u_t(\omega,x)=\tan v_t(\omega,x)
$$
so that $u$ is a renormalized solution of (\ref{Eq11}) in the sense of Definition \ref{Def3}.

(Stability) It follows from the stability in Theorem \ref{Pro3}.

\section{Stochastic Flows with BV Drifts and Constant Diffusions}

Consider the following SDE:
\be
\dif X_t(x)=b(X_t(x))\dif t+\sigma\dif W_t,\ \ X_0=x.\label{SE1}
\ee
In this section, we use Theorem \ref{Th1} to prove the following result.
\bt\label{T1}
Assume that $b$ is a BV vector field and satisfies
$$
\frac{b(x)}{1+|x|},\ \div b(x)\in L^\infty(\mR^d),\ \ b\in \mathrm{BV}_{loc}.
$$
Then there exists a unique almost everywhere stochastic invertible flow to SDE (\ref{SE1}) in the sense of
Definition \ref{Def1}.
\et
\begin{proof}
(Existence): Define $b_n$ as in (\ref{BN}). Let $X_{n,s,t}(x)$ solve the following SDE:
\be
X_{n,s,t}(x)=x+\int^t_sb_n(X_{n,s,r}(x))\dif r +\sigma (W_t-W_s),\ \
\forall t\geq s\geq0.\label{Lp3}
\ee
We divide the proof into two steps.

(Step 1): Fix $t>0$. By Proposition \ref{Pro26}, $v^k_{n,s,t}(x):=X^k_{n,s,t}(x)$
solves the following backward stochastic Kolmogorov equation:
$$
\dif v^k_n+\frac{1}{2}\sigma^{il}\sigma^{jl}\p^2_{ij}v^k_n\dif s+b^i_n\p_i v^k_n\dif s
+\sigma^{il} \p_iv^k_n*\dif W_s=0,\ \ v^k_n|_{s=t}=x^k,
$$
and by Proposition \ref{Pro25}, $u^k_{n,t}(x):=[X^{-1}_{n,0,t}(x)]^k$ solves the following equation:
$$
\dif u^k_n=\frac{1}{2}\sigma^{il}\sigma^{jl}\p^2_{ij}u^k_n
-b^i_n\p_i u^k_n\dif t-\sigma^{il}\p_i u^k_n\dif W^l_t,\ \ u^k_n|_{t=0}=x^k,
$$
where $x^k$ is the $k$-th coordinate of spatial variable $x$.

By Theorem \ref{Th1}, let $v^k_{s,t}$ and $u^k_s$ be the unique renormalized solutions
of the following SPDEs in the sense of Definition \ref{Def3}
\ce
&&\dif v^k+\frac{1}{2}\sigma^{il}\sigma^{jl}\p^2_{ij}v^k\dif s+b^i\p_i v^k\dif s
+\sigma^{il} \p_iv^k*\dif W_s=0,\ \ v^k|_{s=t}=x^k,\\
&&\dif u^k=\frac{1}{2}\sigma^{il}\sigma^{jl}\p^2_{ij}u^k
-b^i\p_i u^k\dif s-\sigma^{il}\p_i u^k\dif W^l_s,\ \ u^k|_{s=0}=x^k.
\de
Then by the stability result in Theorem \ref{Th1}, we have for any $p>1$ and $N>0$,
\be
\lim_{n\to \infty}\mE\left(\sup_{s\in[0,t]}\int_{B_N}|\arctan v^k_{n,s,t}
-\arctan v^k_{s,t}|^p\right)=0\label{L2}
\ee
and
\be
\lim_{n\to \infty}\mE\left(\sup_{s\in[0,t]}\int_{B_N}
|\arctan u^k_{n,s}-\arctan u^k_s|^p\right)=0.\label{L3}
\ee
Define
$$
X_t(\omega, x):=v_{0,t}(\omega,x),\ \ Y_t(\omega,x):=u_t(\omega,x)
$$
Below, we want to show that $X_t(x)$ satisfies {\bf (A)}, {\bf (B)} and {\bf (C)}
of Definition \ref{Def1} and $X_t^{-1}(\omega,x)=Y_t(\omega,x)$.

(Step 2):
By (\ref{L2}), we have for any $p>1$ and $N>0$
$$
\lim_{n\to \infty}\mE\left(\int^t_0\!\!\!\int_{B_N}|\arctan v^k_{n,0,s}
-\arctan v^k_{0,s}|^p\right)=0.
$$
Hence, there exists a subsequence still denoted by $n$ such that
for almost all $(s,\omega,x)\in[0,t]\times\Omega\times\mR^d$ and any $k=1,\cdots,d$
$$
\lim_{n\to\infty}\arctan v^k_{n,0,s}(\omega,x)=\arctan v^k_{0,s}(\omega,x),
$$
i.e.,
$$
\lim_{n\to\infty}X_{n,0,s}(\omega,x)=X_{s}(\omega,x),
$$
as well as for ($P\times\sL$)-almost all $(\omega,x)\in\Omega\times\mR^d$,
\be
\lim_{n\to\infty}X_{n,0,t}(\omega,x)=X_{t}(\omega,x).\label{L5}
\ee
Note that by (\ref{Es5}) and (\ref{EE2}), for any $p\geq 1$,
$$
\mE|\det(\nabla X^{-1}_{n,0,t}(x))|^p\leq e^{C_p t\|\div b_n\|_\infty}\leq
e^{C_pt(\|\div b\|_\infty+\|b/(1+|x|)\|_\infty)}.
$$
By Lemma \ref{Le4}, it is easy to see that {\bf (A)} and {\bf (B)}
of Definition \ref{Def1} hold, and for any $N>0$,
$$
\lim_{n\to\infty}\mE\iint_{B_N}|\div b_n(X_{n,0,s})-\div b(X_s)|=0.
$$
Thus, for ($P\times\sL$)-almost all $(\omega,x)\in\Omega\times\mR^d$,
\be
\det(\nabla X_{n,0,t}(\omega,x))&=&\exp\left\{\int^t_0\div b_n(X_{n,0,s}(\omega,x))\dif s\right\}\no\\
&\stackrel{n\to\infty}{\longrightarrow}&\exp\left\{\int^t_0\div b(X_s(\omega,x))\dif s\right\}
=:\rho_t(\omega,x).\label{L4}
\ee
On the other hand, for fixed $t\geq 0$ and $P$-almost all $\omega\in\Omega$,
it holds that for all $\varphi,\psi\in C_c^+(\mR^d)$,
\be
\int \varphi(u_{n,t}(\omega))\cdot\psi=\int \varphi(X^{-1}_{n,0,t}(\omega))\cdot\psi
=\int \varphi\cdot\psi(X_{n,0,t}(\omega))\cdot\det(\nabla X_{n,0,t}(\omega)).\label{Lp4}
\ee
If necessary, by extracting a subsequence and then taking limits $n\to\infty$ for both sides
of (\ref{Lp4}), by (\ref{L3}), (\ref{L5}) and (\ref{L4}),  we obtain
that for $P$-almost all $\omega\in\Omega$ and all $\varphi,\psi\in C_c^+(\mR^d)$,
$$
\int \varphi(Y_{t}(\omega))\cdot\psi
=\int \varphi\cdot\psi(X_{t}(\omega))\cdot\rho_{t}(\omega).
$$
Thus, by (ii) of Lemma \ref{Le3}, one sees that {\bf (C)} of Definition \ref{Def1} holds.

(Uniqueness): It follows from Propositions \ref{Pro1} and \ref{Pro3}.
\end{proof}
%\br
%It should be noted that Theorem \ref{T1} can be easily obtained by transferring SDE (\ref{SE1}) to ODE. Here, our proof indicates that we can also go back to SDE from stochastic transport equation like from deterministic transport equation to ODE (cf. \cite{Di-Li}).
%\er

\section{Stochastic Flows with Sobolev Drifts and Non-Constant Diffusions}

We first prove the following key estimate.
\bl\label{Le1}
Let $X_t(x)$ and $\hat X_t(x)$ be two almost everywhere stochastic flows of (\ref{SSDE})
corresponding to $(b,\sigma)$ and $(\hat b,\hat \sigma)$ in the sense of
Definition \ref{Def1}, where
$$
b, \hat b\in L^1_{loc}(\mR^d),\ \ |\nabla\hat b|\in (L\log L)_{loc}(\mR^d)
$$
and
$$
\sigma,\hat \sigma\in L^2_{loc}(\mR^d),\ \ |\nabla \hat\sigma|\in L^{2}_{loc}(\mR^d).
$$
Then, for any $T,N,R>0$, there exist constants $C_1, C_2$ given below
such that for all $\delta>0$,
\ce
&&\mE\int_{B_N\cap G^R_T}
\log\left(\frac{\sup_{t\in[0,T]}|X_t-\hat X_t|^{2}}{\delta^2}+1\right)\leq\\
&&\qquad\leq C_1+\frac{C_2}{\delta}
\left(\int_{B_R}|b-\hat b|+\left[\int_{B_R}
|\sigma-\hat \sigma|^2\right]^{1/2}\right),
\de
where
\ce
&&G_T^R(\omega):=\Big\{x\in\mR^d: \sup_{t\in[0,T]}|X_t(\omega,x)|\vee|\hat X_t(\omega,x)|\leq R\Big\},\\
&&C_1:=C_{d,R,N}\cdot T\cdot (K_{T,b,\sigma}+ K_{T,\hat b,\hat\sigma})
\left(1+\int_{B_{2R}}|\nabla\hat b|\log(|\nabla \hat b|+1)+
\left[\int_{B_{2R}}|\nabla\hat \sigma|^{2}
\right]^{\frac{1}{2}}\right),
\de
and $C_2:=C_N\cdot T\cdot K_{T,b,\sigma}$.
Here, $K_{T,b,\sigma}$ is from (\ref{Den}), $C_{d,R,N}$ only depends on $d,R,N$,  and
$C_N$ only depends on $N$.
\el
\begin{proof}
Set
$$
Z_t(x):=X_t(x)-\hat X_t(x).
$$
By It\^o's formula, we have
\ce
\log\left(\frac{|Z_t|^2}{\delta^2}+1\right)
&=&2\int^t_0\frac{\<Z, b(X)-\hat b(\hat X)\>}{|Z|^2+\delta^2}\dif s
+2\int^t_0\frac{\<Z, (\sigma(X)-\hat \sigma(\hat X))\dif W_s\>}{|Z|^2+\delta^2}\\
&&+\int^t_0\frac{\|\sigma(X)-\hat \sigma(\hat X)\|^2}{|Z|^2+\delta^2}\dif s
-2\int^t_0\frac{|(\sigma(X)-\hat \sigma(\hat X))^\mathrm{t}\cdot Z|^2}{(|Z|^2+\delta^2)^2}\dif s\\
&=:&I_1(t)+I_2(t)+I_3(t)+I_4(t).
\de
For $I_1(t)$, we have
$$
I_1(t)\leq\frac{1}{\delta}\int^t_0|b(X)-\hat b(X)|\dif s
+2\int^t_0\frac{|\hat b(X)-\hat b(\hat X)|}{\sqrt{|Z|^2+\delta^2}}\dif s=:I_{11}(t)+I_{12}(t).
$$
Below, we write for a continuous function $f:\mR_+\to\mR$,
$$
f^*(T):=\sup_{t\in[0,T]}|f(t)|.
$$
Noting that
$$
G_T^R(\omega)\subset \{x: |X_t(\omega,x)|\leq R\}\cap\{x: |\hat X_t(\omega,x)|\leq R\},\ \ \forall t\in[0,T],
$$
by (\ref{Den}), we have
$$
\mE\int_{G_T^R}|I^*_{11}(T)|\leq
\frac{1}{\delta}\mE\int^T_0\!\!\!\int_{\{|X|\leq R\}}|b(X)-\hat b(X)|
\leq\frac{\tilde K_{T,b,\sigma}}{\delta}\int_{B_R}|b-\hat b|,
$$
where $\tilde K_{T,b,\sigma}:=T\cdot K_{T,b,\sigma}$,
and by $\sL\circ X\ll\sL$ and $\sL\circ \hat X\ll\sL$,
\ce
\mE\int_{G_T^R}|I^*_{12}(T)|&\stackrel{(\ref{Es2})}{\leq}& C_d\mE\int^T_0\!\!\!\int_{G_T^R}\Big([M_R|\nabla\hat b|](X)
+[M_R|\nabla\hat b|](\hat X)\Big)\\
&\leq& C_d\mE\int^T_0\left(\int_{\{|X|\leq R\}}[M_R|\nabla\hat b|](X)
+\int_{\{|\hat X|\leq R\}}[M_R|\nabla\hat b|](\hat X)\right)\\
&\leq& C_d\cdot(\tilde K_{T,b,\sigma}+\tilde K_{T,\hat b,\hat\sigma}) \int_{B_R}M_R|\nabla\hat b|\\
&\stackrel{(\ref{Es3})}{\leq}& C_{d,R}\cdot (\tilde  K_{T,b,\sigma}
+\tilde K_{T,\hat b,\hat\sigma})
\left(1+\int_{B_{2R}}|\nabla\hat b|\log(|\nabla \hat b|+1)\right).
\de
Hence,
\ce
\mE\int_{G_T^R}|I^*_1(T)|
&\leq& C_{d,R}\cdot (\tilde K_{T,b,\sigma}+\tilde K_{T,\hat b,\hat\sigma})
\left(1+\int_{B_{2R}}|\nabla\hat b|\log(|\nabla \hat b|+1)\right)\\
&&+\frac{\tilde K_{T,b,\sigma}}{\delta}\int^T_0\!\!\!\int_{B_R}|b-\hat b|.
\de
For $I_2(t)$, set
$$
\tau_R(\omega,x):=\inf\Big\{t\geq 0: |X_t(\omega,x)|\vee \hat X_t(\omega,x)>R\Big\},
$$
then
$$
G_T^R(\omega)=\{x:\tau_R(\omega,x)>T\}.
$$
By BDG's inequality, we have
\ce
\mE\int_{B_N\cap G_T^R}|I^*_{2}(T)|&\leq&
\int_{B_N}\mE\left(\sup_{t\in[0,T\wedge\tau_R]}\left|\int^{t}_0\frac{\<Z,(\sigma(X)
-\hat \sigma(\hat X))\dif W_s\>}{|Z|^2+\delta^2}\right|\right)\\
&\leq&C\int_{B_N}\mE\left[\int^{T\wedge\tau_R}_0\frac{|\sigma(X)-\hat \sigma(\hat X)|^2}
{|Z|^2+\delta^2}\dif s\right]^{\frac{1}{2}}\\
&\leq&C_N\left[\mE\int^T_0\!\!\!\int_{\{x:\tau_R(x)\geq s\}}
\frac{|\sigma(X)-\hat \sigma(\hat X)|^2}{|Z|^2+\delta^2}\right]^{\frac{1}{2}}.
\de
As the treatment of $I_1(t)$, we can prove that
\ce
\mE\int_{B_N\cap G_T^R}|I^*_{2}(T)|
\leq C_{d}\cdot (\tilde K_{T,b,\sigma}+\tilde K_{T,\hat b,\hat\sigma})
\left[\int_{B_{2R}}|\nabla\hat \sigma|^{2}\right]^{\frac{1}{2}}
+\frac{C_{N}\cdot\tilde K_{T,b,\sigma}}{\delta}\left[
\int_{B_R}|\sigma-\hat \sigma|^2\right]^{\frac{1}{2}}.
\de
$I_3(t)$ is dealt with similarly and $I_4(t)$ is negative and abandoned.
The proof is thus complete.
\end{proof}
\bl\label{Le2}
Let $\Phi(\omega,x):=\sup_{t\in[0,T]}|X_t(\omega,x)-\hat X_t(\omega,x)|^{2}$. Assume that for some $M>0$,
$$
\int_{B_N\cap G_T^R(\omega)}\log\left(\frac{|\Phi(\omega)|}
{\delta^2}+1\right)\leq M,
$$
where $G^R_T(\omega)$ is as in Lemma \ref{Le1}. Then,
$$
\int_{B_N\cap G_T^R(\omega)}|\Phi(\omega)|\leq \frac{4R^2}{M}+\delta^2(e^{M^2}-1)|B_N|,
$$
where $|B_N|$ denotes the volume of the ball $B_N$.
\el
\begin{proof}
It follows from
$$
\log\left(\frac{|\Phi(\omega,x)|}
{\delta^2}+1\right)\leq M^2\Longrightarrow
|\Phi(\omega,x)|\leq \delta^2(e^{M^2}-1)
$$
and Chebyshev's inequality.
\end{proof}

We introduce the following assumptions on  $b$ and $\sigma$:
\begin{enumerate}[{\bf (H1)}]
\item $b\in L^1_{loc}(\mR^d), |\nabla b|\in (L\log L)_{loc}(\mR^d)$ and $\sigma\in L^2_{loc}(\mR^d),
|\nabla \sigma|\in L^2_{loc}(\mR^d)$.
\item
There exist $b_n,\sigma_n\in C^\infty_b(\mR^d)$ such that
\begin{enumerate}[(i)]
\item For any $R>0$
\be
\lim_{n\to\infty}\int_{B_R}|b_n-b|=0,\ \ \lim_{n\to\infty}\int_{B_R}|\sigma_n-\sigma|^2=0\label{En1}
\ee
and
\be
\sup_n\left(\int_{B_R}|\nabla b_n|(\log(|\nabla b_n|+1))+\int_{B_R}|\nabla\sigma_n|^2\right)<+\infty.\label{En2}
\ee
\item For some $C_1,C_2>0$ independent of $n$,
\be
\|[-\div b_n+\frac{1}{2}\p_i\sigma^{jl}_n \p_j\sigma^{il}_n
+\sigma_n^{il}\p_{ij}^2\sigma^{jl}_n+|\div\sigma_n|^2]^+\|_\infty\leq C_1\label{En3}
\ee
and
\be
\<x,b_n(x)\>_{\mR^d}+2\|\sigma_n(x)\|_{H.S.}^2\leq C_2(|x|^2+1), \ \ \forall x\in\mR^d.\label{En4}
\ee
\end{enumerate}
\end{enumerate}

We are now in a position to prove our main result of this section.
\bt\label{Th11}
Assume that {\bf (H1)} and {\bf (H2)} hold.
Then there exists a unique almost everywhere stochastic
flow of (\ref{SSDE}) in the sense of Definition \ref{Def1}. Moreover, the constant
$K_{T,b,\sigma}$ in (\ref{Den}) is less than $e^{C_1 T}$, where $C_1$ is from (\ref{En3}). In particular,
if $C_1=0$, then $K_{T,b,\sigma}\leq 1$.
\et
\begin{proof}
(Existence): Let $b_n$ and $\sigma_n$ be as in {\bf (H2)}. Let $X_n$ solve the following SDE
$$
\dif X_n=b_n(X_n)\dif t+\sigma_n(X_n)\dif W_t,\ \ X_n|_{t=0}=x.
$$
We want to prove that for any $T, N>0$ and $q\in[1,2)$,
\be
\lim_{n,m\to\infty}\mE\int_{B_N}\sup_{t\in[0,T]}|X_{n,t}(x)-X_{m,t}(x)|^q\dif x=0.\label{C2}
\ee
First of all, by (\ref{En4}), it is standard to prove that
\be
\sup_n\sup_{x\in B_N}\mE\left(\sup_{t\in[0,T]}|X_{n,t}(x)|^2\right)<+\infty.\label{C1}
\ee
Thus, for proving (\ref{C2}), it suffices to prove that for any $\eta>0$,
\be
\lim_{n,m\to\infty}P\left\{\omega: \int_{B_N}
\sup_{t\in[0,T]}|X_{n,t}(\omega,x)-X_{m,t}(\omega,x)|^2\dif x\geq 2\eta\right\}=0.\label{C0}
\ee

Fix $\eps,\eta, T>0$ below and set
$$
\Phi_{n,m}(\omega,x):=\sup_{t\in[0,T]}|X_{n,t}(\omega,x)-X_{m,t}(\omega,x)|^2
$$
and
$$
G^R_{n,m}(\omega):=\Big\{x\in\mR^d: \sup_{t\in[0,T]}
|X_{n,t}(\omega,x)|\vee|X_{m,t}(\omega,x)|\leq R\Big\}.
$$
Then,
\be
&&P\left\{\omega: \int_{B_N}
\Phi_{n,m}(\omega)\geq 2\eta\right\}\leq P\left\{\omega: \int_{B_N\cap G^{R}_{n,m}(\omega)^c}
\Phi_{n,m}(\omega)\geq\eta\right\}\no\\
&&\qquad\qquad+P\left\{\omega: \int_{B_N\cap G^R_{n,m}(\omega)}
\Phi_{n,m}(\omega)\geq\eta\right\}=:I_{n,m}^R+J_{n,m}^R.\label{C3}
\ee
For $I_{n,m}^R$, by Chebyshev's inequality and (\ref{C1}), we may choose $R>0$ large
enough such that for all $n,m\in\mN$,
\be
I_{n,m}^R\leq \frac{1}{\eta}\mE\int_{B_N\cap (G^{R}_{n,m})^c}\Phi_{n,m}
\leq \frac{1}{\eta}\int_{B_N}\Big(\mE\Phi_{n,m}^2 \cdot
P\{\omega: x\notin G^R_{n,m}(\omega)\}\Big)^{\frac{1}{2}}
\leq\eps.\label{C4}
\ee
Fixing such a $R$, we look at $J_{n,m}^R$. Set
$$
\xi^\delta_{n,m}:=\int_{B_N\cap G^R_{n,m}}
\log\left(\frac{\Phi_{n,m}}{\delta^2}+1\right).
$$
By (\ref{Es5}) and (\ref{En3}), we have
\be
\sup_n\sup_{t\in[0,T], x\in\mR^d}\mE|\det(\nabla X^{-1}_{n,t})(x)|^2\leq e^{2TC_1},\label{Es6}
\ee
which yields that the constant $K_{T,b_n,\sigma_n}$ in (\ref{Den}) is bounded by $e^{C_1 T}$. Hence, in Lemma \ref{Le1}, if we choose
$$
\delta=\delta_{n,m}=\int_{B_R}|b_n-b_m|+
\left[\int_{B_R}|\sigma_n-\sigma_m|^2\right]^{\frac{1}{2}},
$$
then by (\ref{En2}), we have for some $C_{T,R,N}$ independent of $n,m$,
$$
\mE \xi^{\delta_{n,m}}_{n,m}\leq C_{T,R,N}.
$$
Thus, there exists an $M_1>0$ such that for all $M\geq M_1$ and all $n,m$,
$$
P(\xi^{\delta_{n,m}}_{n,m}>M)\leq \eps.
$$
Now, by Lemma \ref{Le2} and (\ref{En1}), we may choose $M>M_1\vee 8R^2/\eta$ and $n,m$
large enough such that
$$
\delta_{n,m}<\sqrt{\frac{\eta}{4(e^{M^2}-1)|B_N|}},
$$
which leads to
$$
\Omega^M_{n,m}:=\left\{\omega: \int_{B_N\cap G^R_{n,m}(\omega)}
\Phi_{n,m}(\omega)\geq\eta; \xi^{\delta_{n,m}}_{n,m}(\omega)\leq M\right\}=\emptyset.
$$
Hence, first letting $M$ large enough and then $n,m$ large enough, we obtain
\be
J_{n,m}^R\leq P(\Omega^M_{n,m})+P(\xi^{\delta_{n,m}}_{n,m}>M)\leq\eps.\label{C5}
\ee
Combining (\ref{C3}), (\ref{C4}) and (\ref{C5}), by the arbitrariness of $\eps$,
we get (\ref{C0}) as well as (\ref{C2}).
So, for $q\in(1,2)$,
there exists a stochastic field $X\in L^q_{loc}(\mR^d;L^q(\Omega; C([0,T])))$
such that for any $N>0$
$$
\lim_{n\to\infty}\mE\int_{B_N}\sup_{t\in[0,T]}|X_{n,t}(x)-X_t(x)|^q\dif x=0.
$$
In particular, there is a subsequence still denoted by $n$ such that
for $(P\times\sL)$-almost all $(\omega,x)\in\Omega\times\mR^d$
\be
\lim_{n\to\infty}\sup_{t\in[0,T]}|X_{n,t}(\omega,x)-X_t(\omega,x)|=0.\label{C6}
\ee
In view of (\ref{C1}), (\ref{Es6}) and (\ref{C6}), by Lemma \ref{Le4} and (\ref{En1}),
it is easy to check that $X_t(\omega,x)$ satisfies {\bf (A)} and {\bf (B)}
of Definition \ref{Def1}.

(Uniqueness): Let $X_t(x)$ and $\hat X_t(x)$
be two almost everywhere stochastic flows of (\ref{SSDE}). Then, by Lemma \ref{Le1}, we have
for any $T,N,R>0$ and $\delta>0$,
$$
\mE\int_{B_N\cap G_T^R}\log\left(\frac{\sup_{t\in[0,T]}|X_t-\hat X_t|^2}{\delta^2}+1\right)
\leq C_{T,N,R},
$$
where $C_{T,N,R}$ is independent of $\delta$. Letting $\delta$ go to zero, we obtain
$$
1_{G_T^R(\omega)}(x)\cdot\sup_{t\in[0,T]}|X_t(\omega,x)-\hat X_t(\omega, x)|=0\ \
\mbox{a.e. on $\Omega\times B_N$}
$$
The uniqueness then follows by letting $R\to\infty$.
\end{proof}

The following example is inspired by \cite{Ott, Le-Li1}.

{\bf Example:} Let $d\geq 3$. Consider the following SDE in $\mR^d$ with discontinuous
and degenerate coefficients:
$$
\dif X_t=\frac{\beta X_t}{|X_t|^2}\dif t+\frac{X_t\otimes X_t}{|X_t|^2}\dif W_t,\ \ X_0=x,
$$
where $\beta\geq (4d^2+5d)/(d-2)$.
Define
$$
b(x):=\frac{\beta x}{|x|^2},\ \ \ \
\sigma(x):=\frac{x\otimes x}{|x|^2}
$$
and
$$
b_n(x):=\frac{\beta x}{|x|^2+1/n}, \ \ \sigma_n(x):=\frac{x\otimes x}{|x|^2+1/n}.
$$
By virtue of $d\geq 3$, one sees that for any $q\in(1,3/2)$
$$
|\nabla b|\in L^q_{loc}(\mR^d)\subset(L\log L)_{loc}(\mR^d),\ \
|\nabla\sigma|\in L^2_{loc}(\mR^d).
$$
Thus, {\bf (H1)} is true for $b$ and $\sigma$.

Let us verify {\bf (H2)}. First of all,  (\ref{En1}), (\ref{En2}) and (\ref{En4}) are easily checked.
We look at (\ref{En3}). Noting that
$$
\p_i\sigma^{jl}_n(x)=\frac{\p_i(x^jx^l)(|x|^2+1/n)-2x^ix^jx^l}{(|x|^2+1/n)^2},
$$
we have
$$
\div \sigma^{\cdot l}=\p_i\sigma^{il}(x)=\frac{((d-1)|x|^2+(d+1)/n)x^l}{(|x|^2+1/n)^2}.
$$
Hence,
\ce
\sum_l|\div \sigma^{\cdot l}_n(x)|^2&=&\frac{((d-1)|x|^2+(d+1)/n)^2|x|^2}{(|x|^2+1/n)^4}\\
&\leq&\frac{((d-1)|x|^2+(d+1)/n)^2}{(|x|^2+1/n)^3}\leq\frac{4d^2}{|x|^2+1/n}
\de
and
$$
\p_i\sigma^{jl}_n(x)\p_j\sigma^{il}_n(x)=\frac{(d+3)|x|^2(|x|^2+1/n)^2-8|x|^4/n
-4|x|^6}{(|x|^2+1/n)^{4}}\leq\frac{d+3}{|x|^2+1/n}.
$$
Similarly, we have
$$
\sigma^{il}_n(x)\p_{i}\div\sigma^{\cdot l}_n(x)
=\frac{3(d-1)|x|^4+(d+1)|x|^2/n}{(|x|^2+1/n)^3}
-\frac{4|x|^4((d-1)|x|^2+(d+1)/n)}{(|x|^2+1/n)^4}\leq\frac{4d-2}{|x|^2+1/n}.
$$
Moreover,
$$
\div b_n(x)=\frac{\beta(d-2)}{|x|^2+1/n}+\frac{2\beta}{n(|x|^2+1/n)^2},
$$
Thus, combining the above calculations and by $\beta\geq (4d^2+5d)/(d-2)$, we have
$$
-\div b_n+\frac{1}{2}\p_i\sigma^{jl}_n \p_j\sigma^{il}_n
+\sigma_n^{il}\p_{ij}^2\sigma^{jl}_n+|\div\sigma_n|^2 \leq 0,
$$
and so, (\ref{En3}) holds. Thus, {\bf (H2)} is also true.

We now give two corollaries of Theorem \ref{Th11}.
\bc\label{Cor}
Assume that {\bf (H1)} and {\bf (H2)} hold. Let $Y_0\in L^2(\Omega,\cF_0)$
 be such that $P\circ Y_0\ll\sL$ and the density $\gamma_0\in L^\infty(\mR^d)$.
Then there exists a unique continuous ($\cF_t$)-adapted process $Y_t(\omega)$ such that
\be
\mbox{$P\circ Y_t\ll\sL$ with the density $\gamma_t\in L^\infty_{loc}(\mR_+; L^\infty(\mR^d))$}
\label{Lp7}
\ee
and
$Y_t$ solves
\be
Y_t=Y_0+\int^t_0 b(Y_s)\dif s+\int^t_0\sigma(Y_s)\dif W_s,\ \ \forall t\geq 0.\label{SE0}
\ee
Moreover,
$$
Y_t(\omega)=X_t(\omega,Y_0(\omega)),
$$
where $X_t(x)$ is the unique almost everywhere stochastic flow given in Theorem \ref{Th11}.
\ec
\begin{proof}
As in the proof in the appendix, we can check that $Y_t(\omega):=X_t(\omega,Y_0(\omega))$
solves equation (\ref{SE0}). Moreover, since $X_t(x)$ is independent of $Y_0$,
by (\ref{Den}), we have for any $\varphi\in\cL^+(\mR^d)$ and $t\in[0,T]$,
\ce
\mE\varphi(Y_t)&=&\mE(\mE\varphi(X_t(x))|_{x=Y_0})
=\int_{\mR^d}\mE\varphi(X_t(x))\gamma_0(x)\dif x\\
&\leq&\|\gamma_0\|_\infty\int_{\mR^d}\mE\varphi(X_t(x))\dif x
\leq \|\gamma_0\|_\infty\cdot K_{T,b,\sigma}\int_{\mR^d}\varphi(x)\dif x,
\de
which implies that $P\circ Y_t\ll\sL$ and the density $\gamma_t$ satisfies
$$
\sup_{t\in[0,T]}\|\gamma_t\|_\infty\leq \|\gamma_0\|_\infty\cdot K_{T,b,\sigma}.
$$

Let us now look at the uniqueness. Let $\hat Y_t$ be another solution of (\ref{SE0}) with
$\hat Y_0=Y_0$ and satisfy that
\be
P\circ \hat Y_t\ll\sL\ \ \mbox{with the density $\hat \gamma_t\in
L^\infty_{loc}(\mR_+; L^\infty(\mR^d))$}.\label{Lp8}
\ee
It is now standard to prove that for any $T>0$,
\be
\mE\left(\sup_{t\in[0,T]}|Y_t|^2\right)+\mE\left(\sup_{t\in[0,T]}|\hat Y_t|^2\right)<+\infty.\label{Po1}
\ee
Set
$$
Z_t:=Y_t-\hat Y_t
$$
and for $R>0$
$$
\tau_R:=\inf\{t\geq 0: |Y_t|\vee|\hat Y_t|\geq R\}.
$$
Then by (\ref{Po1}), we have
$$
P\left\{\omega:\lim_{R\to\infty}\tau_R(\omega)=+\infty\right\}=1.
$$
As in the proof of Lemma \ref{Le1}, we have
\be
\mE\log\left(\frac{|Z_{t\wedge\tau_R}|^2}{\delta^2}+1\right)
&\leq&2\mE\int^{t\wedge\tau_R}_0\frac{\<Z, b(Y)-b(\hat Y)\>}{|Z|^2+\delta^2}\dif s
+\mE\int^{t\wedge\tau_R}_0\frac{\|\sigma(Y)-\sigma(\hat Y)\|^2}{|Z|^2+\delta^2}\dif s\no\\
&\stackrel{(\ref{Es2})}{\leq}&C\mE\int^{t\wedge\tau_R}_0([M_R|\nabla b|](Y)
+[M_R|\nabla b|](\hat Y))\dif s+C_T\label{Ep4}\\
&\leq&C\int^{t}_0(\mE (1_{|Y|\leq R}\cdot[M_R|\nabla b|](Y))+\mE(1_{|\hat Y|\leq R}
\cdot[M_R|\nabla b|](\hat Y)))\dif s+C_T\no\\
&\stackrel{(\ref{Lp7}) (\ref{Lp8})}{\leq}&C_T\int_{|y|\leq R}[M_R|\nabla b|](y)\dif y+C_T\no\\
&\stackrel{(\ref{Es3})}{\leq}&C_T\int_{|y|\leq R}|\nabla b|(y)\log(|\nabla b(y)|+1)\dif y+C_T,\no
\ee
which yields the uniqueness by first letting $\delta\to 0$ and then $R\to\infty$.
\end{proof}
\bc\label{Cor1}
In addition to {\bf (H1)} and {\bf (H2)}, we also assume that for some $q>d$,
$$
|\nabla b|\in L^q_{loc}(\mR^d).
$$
Let $Y_0\in L^2(\Omega,\cF_0)$ be such that $P\circ Y_0\ll\sL$ and the
density $\gamma_0\in L^\infty(\mR^d)$.
Then  $Y_t(\omega):=X_t(\omega, Y_0(\omega))$ uniquely solves SDE (\ref{SE0}),
where $X_t(x)$ is the unique almost everywhere stochastic flow given in Theorem \ref{Th11}.
\ec
\begin{proof}
Following the proof of Corollary \ref{Cor}, we only need to prove the uniqueness.
Let $\hat Y$ be another solution of SDE (\ref{SE0}) with the same initial value $\hat Y_0=Y_0$.
Choosing $q'\in(d,q)$, and using (\ref{Es02}) in (\ref{Ep4}), we have
\ce
\mE\log\left(\frac{|Z_{t\wedge\tau_R}|^2}{\delta^2}+1\right)
&\leq&C_{q'}\mE\int^{t\wedge\tau_R}_0[M_R|\nabla b|^{q'}]^{1/q'}(Y)
\dif s+C_T\\
&\stackrel{(\ref{Lp7})}{\leq}&C_{q',T}\int_{|y|\leq R}[M_R|\nabla b|^{q'}]^{1/q'}(y)\dif y+C_T\no\\
&\stackrel{(\ref{Es30})}{\leq}&C_{q',T,q,R}\int_{|y|\leq R}|\nabla b|^q(y)\dif y+C_T,
\de
which in turn implies the uniqueness as Corollary \ref{Cor}.
\end{proof}

\section{Proofs of Main Results}

We first give:

{\it Proof of Theorem \ref{Main2}:}
Under (\ref{BB}) and (\ref{BB0}), it has been proven in Theorem \ref{T1}.
We now consider the case of (\ref{BB}) and (\ref{SI}).
Let us define $b_n:=b*\varrho_n\cdot\chi_n$ and $\sigma_n:=\sigma*\varrho_n\cdot\chi_n$
as in (\ref{BN}). Note that as in estimating (\ref{EE2}),
\ce
|\nabla b_n|&\leq&|\nabla b|*\varrho_n\cdot\chi_n+|b|*\varrho_n\cdot|\nabla\chi_n|\\
&\leq&|\nabla b|*\varrho_n+2\|\nabla\chi\|_\infty\cdot\|b/(1+|x|)\|_\infty\\
&=:&|\nabla b|*\varrho_n+C_1.
\de
If we define
$$
\Psi(r):=(r+C_1)\log(r+C_1+1),
$$
then $r\to\Psi(r)$ is a convex function on $\mR_+$.
Thus, by Jensen's inequality, we have for any $R>0$,
\be
\int_{B_R}|\nabla b_n|\log(|\nabla b_n|+1)
\leq\int_{B_R}\Psi(|\nabla b|*\varrho_n)\leq\int_{B_R}\Psi(|\nabla b|)*\varrho_n
\leq\int_{B_R}\Psi(|\nabla b|).\label{PL2}
\ee
Moreover, by (\ref{BB}) and (\ref{SI}), it is easy to check that
\be
\sup_n\left(\big\|\frac{|b_n|}{1+|x|}\big\|_\infty+\|\div b_n\|_\infty+\|\nabla\sigma_n\|_\infty+
\||\sigma_n|\cdot|\nabla\div\sigma_n|\|_\infty\right)<+\infty.\label{PL1}
\ee
Hence, {\bf (H1)} and {\bf (H2)} hold.

By Theorem \ref{Th11}, there exists a unique
almost everywhere stochastic flow. Following the proof of Theorem \ref{Th11},
we only need to check {\bf (C)} of Definition \ref{Def1}.

Fix a $T>0$ and let
$$
\rho_{n}:=\exp\left\{\int^T_0\Big(\div b_n
-\frac{1}{2}\p_i\sigma^{jl}_n\p_j\sigma^{il}_n\Big)(X_n)\dif s
+\int^T_0\div \sigma_n(X_n)\dif W_s\right\}.
$$
As in Lemma \ref{Le5} and by (\ref{PL1}), we have for any $p\geq 1$,
\be
\sup_{n\in\mN}\sup_{x\in\mR^d}\mE|\rho_n(x)|^p<+\infty.\label{PP3}
\ee
In view of (\ref{C1}), (\ref{Es6}) and (\ref{C6}), by Lemma \ref{Le4}, we have for any $N>0$,
\ce
&&\lim_{n\to\infty}\mE\int^T_0\!\!\!\int_{B_N}|\div b_n(X_{n})-\div b(X)|=0,\\
&&\lim_{n\to\infty}\mE\int^T_0\!\!\!\int_{B_N}|\p_i\sigma^{jl}_n\p_j\sigma^{il}_n(X_{n})
-\p_i\sigma^{jl}\p_j\sigma^{il}(X)|=0,\\
&&\lim_{n\to\infty}\mE\int_{B_N}
\left|\int^T_0(\div\sigma_n(X_{n})-\div \sigma(X))\dif W_s\right|=0.
\de
So, there is a subsequence still denoted by $n$ such that for almost all $(\omega,x)$,
\be
\lim_{n\to\infty}\rho_{n}(\omega,x)=\rho_T(\omega,x),\label{PP4}
\ee
where $\rho_T(x)$ is defined by (\ref{Rho}). By (\ref{PP3}) and (\ref{PP4}), we further have
for any $p\geq 1$ and $N>0$,
\be
\lim_{n\to\infty}\mE\int_{B_N}|\rho_{n}-\rho_T|^p=0.\label{U1}
\ee

Now, let $Y_n$ solve the following SDE
$$
\dif Y_n=-\hat b_n(Y_n)\dif t+\sigma_n(Y_n)\dif W^T_t,\ \ Y_n|_{t=0}=x,
$$
where $\hat b^i_n=b^i_n-\sigma^{jl}_n\p_j\sigma^{il}_n$ and $W^T_t:=W_{T-t}-W_T$.
As in the proof of Theorem \ref{Th11}, there exists
$$
Y\in L^2_{loc}(\mR^d; L^2(\Omega; C([0,T])))
$$
such that for any $N>0$,
\be
\lim_{n\to\infty}\mE\int_{B_N}\sup_{t\in[0,T]}|Y_{n,t}(x)-Y_t(x)|^2\dif x=0.\label{U2}
\ee
Note that for any $\varphi,\psi\in C_c^+(\mR^d)$ (see the proof of Lemma \ref{Le5}),
\be
\int \varphi(Y_{n,T}(\omega))\cdot\psi=\int \varphi\cdot\psi(X_{n,T}(\omega))\cdot\rho_{n}(\omega),
\ \ P-a.s.\label{L01}
\ee
By (\ref{C6}), (\ref{U1}) and (\ref{U2}),  if necessary, extracting a subsequence and
then taking limits $n\to\infty$ in $L^1(\Omega)$ for both sides of (\ref{L01}), we get
that for all $\varphi,\psi\in\sC\subset C^+_c(\mR^d)$ and $P$-almost all $\omega\in\Omega$,
\be
\int \varphi(Y_{T}(\omega))\cdot\psi=\int \varphi\cdot\psi(X_{T}(\omega))\cdot\rho_{T}(\omega).\label{L1}
\ee
Since $\sC$ is countable, one may find a common null set $\Omega'\subset\Omega$
such that (\ref{L1}) holds for all $\omega\notin\Omega'$ and $\varphi,\psi\in\sC$.
Thus, by (ii) of Lemma \ref{Le3}, one sees that {\bf (C)} of Definition \ref{Def1} holds.

We next give:

{\it Proof of Theorem \ref{Main1}:}
We follow the classical Krylov-Bogoliubov's method.
Let $Y_0$ be an $\cF_0$-measurable $\mR^d$-valued random variable. Suppose that the probability law of
$Y_0$ is absolutely continuous with respect to $\sL$ with the density $\gamma_0\in L^\infty(\mR^d)$.
Define $Y_t(\omega):=X_t(\omega,Y_0(\omega))$ and
$$
\mu_n(\varphi):=\frac{1}{n}\int^n_0\mE\varphi(Y_s)\dif s=
\frac{1}{n}\int^n_0\mE[(\mT_s\varphi)(X_0)]\dif s,
$$
where $\{X_s(x),x\in\mR^d\}_{t\geq 0}$ is the unique almost everywhere stochastic flow of (\ref{SSDE}).

Noting that $Y_t(\omega)$ solves the following SDE (see Corollary \ref{Cor})
$$
Y_t=Y_0+\int^t_0b(Y_s)\dif s+\int^t_0\sigma(Y_s)\dif W_s,
$$
by (\ref{CO})  and It\^o's formula, it is standard to prove that
$$
\mE|Y_t|^2\leq \mE|Y_0|^2\ \mbox{ or }\
\frac{1}{t}\int^t_0\mE|Y_s|^2\dif s\leq\frac{\mE|Y_0|^2}{C_1t}
+\frac{C_2}{C_1}.
$$
From this, we derive that the family of probability measures $\mu_n$ is tight.

On the other hand, for any $\varphi\in\cL^+(\mR^d)$, we have
\ce
\mu_n(\varphi)&=&\frac{1}{n}\int^n_0\!\!\!\int_{\mR^d}\mT_s\varphi(x)\cdot\gamma_0(x)\dif x\dif s\\
&\leq&\|\gamma_0\|_\infty\frac{1}{n}\int^n_0\!\!\!\int_{\mR^d}\mT_s\varphi(x)\dif x\dif s\\
&\stackrel{(\ref{Den})}{\leq}&\|\gamma_0\|_\infty\cdot K_{b,\sigma}\cdot\int_{\mR^d}\varphi(x)\dif x,
\de
which means that
$$
\mu_n\ll\sL
$$
and the density $\gamma_n$ satisfies
$$
\|\gamma_n\|_\infty\leq \|\gamma_0\|_\infty\cdot K_{b,\sigma}.
$$
Hence, there exists a subsequence $n_k$, $\gamma\in L^\infty(\mR^d)$ and a probability measure $\mu$
such that
$$
\gamma_{n_k}\mbox{ weakly $*$ converges to $\gamma$ in $L^\infty(\mR^d)$}
$$
and $\mu_{n_k}$weakly converges to $\mu$ in the sense that for any $\varphi\in C_b(\mR^d)$
$$
\lim_{k\to\infty}\int_{\mR^d}\varphi(x)\mu_{n_k}(\dif x)=\int_{\mR^d}\varphi(x)\mu(\dif x).
$$
Since for all $\varphi\in C_c(\mR^d)$,
$$
\int_{\mR^d}\varphi(x)\mu(\dif x)=\int_{\mR^d}\varphi(x)\gamma(x)\dif x.
$$
we have $\mu(\dif x)=\gamma(x)\dif x$.

Let us verify (\ref{Eq}). For $\varphi\in L^1(\mR^d)$ and $t\geq 0$, since $\mT_t\varphi\in L^1(\mR^d)$,
we have
\ce
\int_{\mR^d} \mT_t\varphi(x)\gamma(x)\dif x&=&\lim_{k\to\infty}
\int_{\mR^d} \mT_t\varphi(x)\gamma_{n_k}(x)\dif x=
\lim_{k\to\infty}\frac{1}{n_k}\int^{n_k}_0\!\!\!\int_{\mR^d}
\mT_s\mT_t\varphi(x)\gamma_0(x)\dif x\dif s\\
&=&\lim_{k\to\infty}\frac{1}{n_k}\int^{n_k}_0\!\!\!\int_{\mR^d}
\mT_{t+s}\varphi(x)\gamma_0(x)\dif x\dif s\\
&=&\lim_{k\to\infty}\frac{1}{n_k}\left(\int^{n_k}_0\!\!\!\int_{\mR^d}
+\int^{n_k+t}_{n_k}\!\!\!\int_{\mR^d}-\int^t_0\!\!\!\int_{\mR^d}\right)\mT_s\varphi(x)\gamma_0(x)\dif x\dif s\\
&=&\lim_{k\to\infty}\int_{\mR^d}\varphi(x)\gamma_{n_k}(x)\dif x=\int_{\mR^d}\varphi(x)\gamma(x)\dif x.
\de
The proof is thus complete.

\section{Appendix}

Before proving Proposition \ref{Pr2}, we need the following simple lemma.
\bl\label{Le6}
Let $\sG$ and $\sA$ be two independent $\sigma$-subalgebras of $\cF$.
Let $G:\Omega\times\mR^d\to\mR$ be a bounded $\sG\times\cB(\mR^d)$-measurable function
and $X:\Omega\times\mR^d\to\mR^d$ a  $\sA\times\cB(\mR^d)$-measurable mapping.
Suppose that for $P$-almost all $\omega$, $\sL\circ X(\omega,\cdot)\ll\sL$. Then
for $\sL$-almost all $x\in\mR^d$,
\be
\mE (G(\cdot,X(\cdot,x))|\sA)=(\mE G(\cdot,y))|_{y=X(\cdot,x)}.\label{LL2}
\ee
\el
\begin{proof}
Define $G_\eps(\omega,y):=G(\omega,\cdot)*\varrho_\eps(y)$, where $\varrho_\eps$ is a
family of regularized kernel functions as in Section 4. It is easy to see that
\be
\mE(G_\eps(\cdot,X(\cdot,x))|\sA)=(\mE G_\eps(\cdot,y))|_{y=X(\cdot,x)}.\label{LL1}
\ee
Since for ($P\times\sL$)-almost all $(\omega,y)\in\Omega\times\mR^d$,
$$
\lim_{\eps\to 0}G_\eps(\omega,y)=G(\omega,y)
$$
and
$$
(P\times\sL)\circ(\cdot,X(\cdot,\cdot))\ll P\times\sL.
$$
By taking limits $\eps\to 0$ for both sides of (\ref{LL1}), we get (\ref{LL2}).
\end{proof}

{\it Proof of Proposition \ref{Pr2}:}
Consider the case of almost everywhere stochastic invertible flow. Fix an $s>0$ below.
By {\bf (B)} of Definition \ref{Def1}, one sees that
\be
(P\times\sL)\circ(\theta_s(\cdot),X_s(\cdot,\cdot))\ll P\times\sL.\label{PP6}
\ee
Therefore, there exists a null set $A_s\subset\Omega\times\mR^d$ such that
for all $(\omega,x)\notin A_s$,
$$
\tilde X_t(\omega,x):=
\left\{\begin{aligned}
&X_t(\omega,x),\ &  t\in[0,s],\ \ \\
&X_{t-s}(\theta_s\omega,X_s(\omega,x)),\ & \ \ t\in[s,\infty)
\end{aligned}
\right.
$$
is well defined.
We now check that $\tilde X$ still satisfies {\bf (A)}, {\bf (B)} and  {\bf (C)}
of Definition \ref{Def1}.

{\bf Verification of (A) for $\tilde X$}: It is clear that for $\sL$-almost all $x\in\mR^d$,
$t\mapsto\tilde X_t(x)$ is  a continuous and ($\cF_t$)-adapted process.
We just need to show that for any $t>s$,
\be
\int^t_s|b(\tilde X_r(x))|\dif r+\int^t_s|\sigma(\tilde X_r(x))|^2\dif r<+\infty,
\ \ (P\times\sL)-a.e.,\label{PP7}
\ee
and for $\sL$-almost all $x\in\mR^d$,
\be
\tilde X_t(x)=X_s(x)+\int^t_sb(\tilde X_r(x))\dif r+\int^t_s\sigma(\tilde X_r(x))
\dif W_s,\ \ P-a.s.\label{PP8}
\ee
First of all, by (\ref{PP6}) it is easy to see that (\ref{PP7}) is true.
We look at (\ref{PP8}). Write
$$
Y_{s,t}(\omega,x):=X_{t-s}(\theta_s\omega,x),\ \ t\geq s
$$
and for $M>0$, set
$$
\tau_M(\omega,x):=\inf\left\{t\geq 0: \int^t_0|\sigma(X_r(\omega,x))|^2\dif r>M\right\}.
$$
Then for $\sL$-almost all $x,y\in\mR^d$,
\be
\mbox{$\tau_M(\theta_s(\cdot),y)$ and $Y_{s,t}(\omega,y)$ are independent of $X_s(x)$.}\label{P0}
\ee
By {\bf (A)} for $ X$ and (\ref{PP6}), we have for ($P\times\sL$)-almost all
$(\omega,x)\in\Omega\times\mR^d$,
\be
\lim_{M\to\infty}\tau_M(\theta_s(\omega),X_s(\omega,x))=+\infty.\label{PP9}
\ee
Observe that $Y_{s,t}(x)$ solves
$$
Y_{s,t}(x)=x+\int^t_sb(Y_{s,r}(x))\dif r+\int^t_s\sigma(Y_{s,r}(x))\dif W_r,\ \ t\geq s.
$$
For verifying (\ref{PP8}), by (\ref{PP9})
it suffices to show that for $\sL$-almost all $x\in\mR^d$,
\be
\int^{t\wedge\tau_M(\theta_s(\cdot),y)}_s\sigma(Y_{s,r}(y))\dif W_r\Big|_{y=X_s(x)}=
\int^{t\wedge\tau_M(\theta_s(\cdot),X_s(x))}_s\sigma(Y_{s,r}(X_s(x)))\dif W_r, \ P-a.s.\label{Ep2}
\ee
We extend $\sigma(Y_{s,r}(y))=0$ for $r<s$ and define for $h>0$
$$
f^h_r(y):=\frac{1}{h}\int^r_{r-h}\sigma(Y_{s,r'}(y))\dif r'.
$$
Then $r\to f^h_r(y)$ is a continuous and ($\cF_t$)-adapted process and
\be
\int^t_s|f^h_r(y)|^2\dif r\leq\int^t_s|\sigma(Y_{s,r}(y))|^2\dif r,\ \
\lim_{h\to 0}\int^t_s|f^h_r(y)-\sigma(Y_{s,r}(y))|^2\dif r=0.\label{P2}
\ee
Hence, for any $R>0$, by (\ref{P0}) and Lemma \ref{Le6}, we have
\be
&&\mE\int_{|X_s(x)|\leq R}\left|\int^{t\wedge\tau_M(\theta_s(\cdot),X_s(x))}_s(f^h_r(X_s(x))
-\sigma(Y_{s,r}(X_s(x))))\dif W_r\right|^2\dif x\no\\
&&\quad=\int\mE\left|\int^{t\wedge\tau_M(\theta_s(\cdot),X_s(x))}_s(f^h_r(X_s(x))
-\sigma(Y_{s,r}(X_s(x))))\cdot 1_{\{|X_s(x)|\leq R\}}\dif W_r\right|^2\dif x\no\\
&&\quad=\int
\mE\left(\int^{t\wedge\tau_M(\theta_s(\cdot),X_s(x))}_s|f^h_r(X_s(x))
-\sigma(Y_{s,r}(X_s(x)))|^2\cdot 1_{\{|X_s(x)|\leq R\}}\dif r\right)\dif x\no\\
&&\quad=\int\mE\left(
\mE\left(\int^{t\wedge\tau_M(\theta_s(\cdot),y)}_s|f^h_r(y)
-\sigma(Y_{s,r}(y))|^2\cdot 1_{\{|y|\leq R\}}\dif r\right)\Bigg|_{y=X_s(x)}\right)\dif x\no\\
&&\quad\stackrel{(\ref{Den})}{\leq} K_{T,b,\sigma}\int_{B_R}\mE
\left(\int^{t\wedge\tau_M(\theta_s(\cdot),y)}_s|f^h_r(y)
-\sigma(Y_{s,r}(y))|^2 \dif r\right)\dif y\no\\
&&\quad\to 0,\ \ \mbox{ as $h\to 0$},\label{P3}
\ee
where the last step is due to (\ref{P2}) and the dominated convergence theorem.
Similarly, we can prove that
$$
\lim_{h\to 0}\mE\int_{|X_s(x)|\leq R}\left(\int^{t\wedge\tau_M(\theta_s(\cdot),y)}_s(f^h_r(y)
-\sigma(Y_{s,r}(y)))\dif W_r\Bigg|_{y=X_s(x)}\right)^2\dif x=0.
$$
Thus, for proving (\ref{Ep2}), we only need to prove that for fixed $h>0$,
\be
\int^{t\wedge\tau_M(\theta_s(\cdot),y)}_sf^h_r(y)\dif W_r\Big|_{y=X_s(x)}=
\int^{t\wedge\tau_M(\theta_s(\cdot),X_s(x))}_sf^h_r(X_s(x))\dif W_r, \ P-a.s.\label{Ep22}
\ee
Let $\Delta_n=\{s=r_0<r_1<\cdots<r_n=t\}$ be a division of $[s,t]$. Write
$$
F^h_n(y):=\sum_{r_k\in\Delta_n\setminus\{r_n\} }f^h_{r_k}(y)(W_{r_{k+1}}-W_{r_k})
\cdot 1_{r_k\leq\tau_M(\theta_s(\cdot),y)}
$$
and
$$
F^h(y):=\int^{t\wedge\tau_M(\theta_s(\cdot),y)}_s f^h_r(y)\dif W_r.
$$
Then $F^h_n(y)$ and $F^h(y)$ are independent of $X_s(x)$ and for $\sL$-almost all $y\in\mR^d$,
\be
\mE|F_n^h(y)|^2\leq C_{h,M},\ \
\lim_{|\Delta_n|\to 0}\mE|F^h_n(y)-F^h(y)|^2=0,\label{PP10}
\ee
where $|\Delta_n|:=\min_{r_k\in\Delta_n\setminus\{r_n\}}|r_{k+1}-r_k|$. Thus, as in estimating (\ref{P3}),
by (\ref{Den}) and (\ref{PP10}), we have
$$
\lim_{|\Delta_n|\to 0}\mE\int_{|X_s(x)|\leq R}
\left|F^h_n(X_s(x))-\int^{t\wedge\tau_M(\theta_s(\cdot),X_s(x))}_s
f^h_r(X_s(x))\dif W_r\right|^2\dif x=0
$$
and
$$
\lim_{|\Delta_n|\to 0}\mE\int_{|X_s(x)|\leq R}\left(F^h_n(X_s(x))-
\int^{t\wedge\tau_M(\theta_s(\cdot),y)}_s
f^h_r(y)\dif W_r\Bigg|_{y=X_s(x)}\right)^2\dif x=0,
$$
which in turn yields (\ref{Ep22}).

{\bf Verification of (B) for $\tilde X$}:
By (\ref{P0}) and Lemma \ref{Le6}, we have for any bounded measurable function $\varphi$,
$$
\mE\varphi(Y_{s,t}(X_s(x)))=\mE\big(\mE\varphi(Y_{s,t}(y))|_{y=X_s(x)}\big).
$$
Hence, by (\ref{Den}),  we have for any $s\leq t\leq T$
$$
\int_{\mR^d}\mE\varphi(\tilde X_t(x))\dif x
=\int_{\mR^d}\mE\varphi(Y_{s,t}(X_s(x)))\dif x\leq K_{T,b,\sigma}
\int_{\mR^d}\mE\varphi(Y_{s,t}(y))\dif y
\leq K_{T,b,\sigma}^2 \int_{\mR^d}\varphi(x)\dif x.
$$

{\bf Verification of (C) for $\tilde X$}: Fixing $t\geq s$, we have for  $\varphi\in\cL^+(\mR^d)$,
\ce
\int_{\mR^d}\varphi(\tilde X^{-1}_s(\omega,x))\dif x
&=&\int_{\mR^d}\varphi(X^{-1}_s(\omega,X^{-1}_{t-s}(\theta_s\omega,x)))\dif x\\
&=&\int_{\mR^d}\varphi(X^{-1}_s(\omega,x))\rho_{t-s}(\theta_s\omega,x)\dif x\\
&=&\int_{\mR^d}\varphi(x)\rho_{t-s}(\theta_s\omega,X_s(\omega,x))
\rho_s(\omega,x)\dif x.
\de
Noticing that
$$
\rho_{t-s}(\theta_s\cdot,x)=
\exp\left\{\int^t_s\Big[\div b
-\frac{1}{2}\p_i\sigma^{jl}\p_j\sigma^{il}\Big](Y_{s,r}(x))\dif r+\int^t_s
\div \sigma(Y_{s,r}(x))\dif W_r\right\},
$$
as in verifying (\ref{PP8}), we have
$$
\rho_{t-s}(\theta_s\cdot,X_s(x))=
\exp\left\{\int^t_s\Big[\div b
-\frac{1}{2}\p_i\sigma^{jl}\p_j\sigma^{il}\Big](\tilde X_r(x))\dif r+\int^t_s
\div \sigma(\tilde X_r(x))\dif W_r\right\}.
$$
Thus,
\ce
&&\qquad\tilde\rho_t(x)=\rho_{t-s}(\theta_s\cdot,X_s(x))\rho_s(x)=\\
&&=\exp\left\{\int^t_0\Big[\div b
-\frac{1}{2}\p_i\sigma^{jl}\p_j\sigma^{il}\Big](\tilde X_s(x))\dif s+\int^t_0
\div \sigma(\tilde X_s(x))\dif W_s\right\}.
\de

\vspace{5mm}

Finally, by the uniqueness, we have for $(P\times\sL)$-almost all $(\omega,x)\in\Omega\times\mR^d$,
$$
\tilde X_t(\omega,x)=X_t(\omega,x),\ \forall t\geq 0,
$$
that is, (\ref{Ep3}) holds.

{\bf Markov Property (\ref{Ep30})}: It follows from (\ref{Ep3})
 and Lemma \ref{Le6} as well as the independence of $X_{t}(\theta_s\cdot,x)$ and $\cF_s$.

\vspace{2mm}

{\bf Acknowledgements:}

The author would like to thank Professor Benjamin Goldys for
providing him an excellent environment to work in the University of New South Wales.
His work is supported by ARC Discovery grant DP0663153 of Australia and
NSF of China (No. 10871215).

\end{document}